\documentclass[fleqn,12pt]{wlscirep}
\usepackage{caption}
\usepackage{amssymb}
\usepackage{hyperref}       % hyperlinks
\usepackage{url}            % simple URL typesetting
\usepackage{amsfonts}       % blackboard math symbols
\usepackage{lineno}
\title{Identifying heterogeneous micromechanical properties of biological tissues via physics-informed neural networks}

\author[1,2]{Wensi Wu} % wuw4@chop.edu
\author[3,4]{Mitchell Daneker} % mitchell.daneker@yale.edu
\author[5]{Kevin T. Turner} % ktturner@seas.upenn.edu
\author[1,2]{Matthew A. Jolley} % jolleym@chop.edu
\author[3,*]{Lu Lu}

\affil[1]{Department of Anesthesiology and Critical Care Medicine, Children's Hospital of Philadelphia, Philadelphia, PA 19104}
\affil[2]{Division of Cardiology, Children's Hospital of Philadelphia, Philadelphia, PA 19104}
\affil[3]{Department of Statistics and Data Science, Yale University, New Haven, CT 06511}
\affil[4]{Department of Chemical and Biochemical Engineering, University of Pennsylvania, Philadelphia, PA 19104}
\affil[5]{Department of Mechanical Engineering and Applied Mechanics, University of Pennsylvania, Philadelphia, PA 19104}
\affil[*]{Corresponding author. Email: lu.lu@yale.edu}


\begin{abstract}
The heterogeneous micromechanical properties of biological tissues have profound implications across diverse medical and engineering domains. However, identifying full-field heterogeneous elastic properties of soft materials using traditional engineering approaches is fundamentally challenging due to difficulties in estimating local stress fields. Recently, there has been a growing interest in using data-driven models to learn full-field mechanical responses such as displacement and strain from experimental or synthetic data. However, research studies on inferring full-field elastic properties of materials, a more challenging problem, are scarce, particularly for large deformation, hyperelastic materials. Here, we propose a physics-informed machine learning approach to identify the elasticity map in nonlinear, large deformation hyperelastic materials. We evaluate the prediction accuracies and computational efficiency of physics-informed neural networks (PINNs) by inferring the heterogeneous elasticity maps across three materials with structural complexity that closely resemble real tissue patterns, such as brain tissue and tricuspid valve tissue. We further applied our improved architecture to three additional examples of breast cancer tissue and extended our analysis to three hyperelastic constitutive models: Neo-Hookean, Mooney Rivlin, and Gent. Our selected network architecture consistently produced highly accurate estimations of heterogeneous elasticity maps, even when there was up to 10\% noise present in the training data.
\end{abstract}

\begin{document}
\flushbottom
\maketitle
\thispagestyle{empty}
\section{Introduction}
Biological tissues exhibit a vast range of structural and mechanical heterogeneity~\cite{Wei2023, Doksoo2023}. These properties have significant medical implications for both normal physiological function and the treatment of diseases. Research studies have shown that tissue biomechanical properties are strongly correlated to cardiovascular diseases, including aortic disscection~\cite{xuan2021}, aortic aneurysm~\cite{Burris2022, Lin2022, Klaas2023}, myocardial infarction~\cite{Chang2023}, and degenerated heart valves~\cite{Sadeghinia2023}. Further, tissue and cell stiffness also play a critical role in identifying cancer invasion~\cite{Micalet2023}, tissue regenerative capacity~\cite{Seifert2012, Maden2018, Notari2018, Harn2021}, host reaction upon implantation~\cite{Paszek2005, Engler2006, Sommakia2014, Wilson2023}, and fibrosis~\cite{Wu2020, Henderson2020}. Computational analysis, such as finite element modeling, provides a valuable non-invasive method to analyze the formation and progression underlying these medical phenomena~\cite{Wu2022, Wu2023_2, Zhang2021, Kong2020}. Yet, translating finite element computational analysis to inform medical decision-making still requires further model validation, which necessitates precise knowledge of the heterogeneous mechanical property distribution in biological tissues. 

There exist several methods for identifying the mechanical properties of materials, including atomic force microscopy (AFM) indentation, inverse finite element analysis, or virtual fields method. AFM is particularly effective for examining heterogeneous structures and properties of specimens at the cellular level~\cite{dufrene2017, Amo2017}. In particular, a cantilever probe scans the tissue samples and applies perpendicular pressure to measure the resulting deformation. The ratio of the applied pressure and deformation is subsequently used to approximate the tissue stiffness at each probed location. Although AFM produces high-resolution maps of elastic moduli, the indentation loading mode does not fully replicate the physiological conditions of living tissues~\cite{Rosalia2023, Guimaraes2020}. Therefore, the mechanical behaviors of biological tissues may be misrepresented.

In inverse finite element analysis, a finite element model is constructed to represent the geometry of the specimen being tested. The model is then subjected to appropriate boundary and loading conditions to simulate those in the experimental setup of a mechanical test. An initial guess of material constitutive parameters is applied to the finite element model and adjusted iteratively until the resulting mechanical responses agree with those measured in experiments. The estimation of mechanical parameters requires values of stress-strain pairs. While it is relatively straightforward to approximate stress distribution on the boundary of a testing specimen, extracting internal stress is non-trivial and often relies on a complex multiscale modeling approach~\cite{Bruno20216, Salvati2017, Rokos2023}. As such, inverse finite element analysis has primarily been confined to deducing the macroscopic homogeneous elastic properties of the material~\cite{deplano2016, Ma2020, Laville2020}. The virtual fields method is another popular inverse method for identifying the full-field mechanical parameters of heterogeneous materials~\cite{pierron2012, Luetkemeyer2021}. The virtual fields method reconstructs the full-field stress data leveraging full-field displacement measurement provided by digital image correlation in conjunction with a preselected virtual displacement/velocity field function. The full-field material parameters are subsequently optimized by satisfying the Principle of Virtual Work equilibrium. It is important to note that the accuracy of parameter predictions is highly dependent on the selection of the virtual field function~\cite{Deng2023}.

Recently, there has been increased interest in applying data-driven methods to discover unknown material model parameters. Researchers have demonstrated high prediction accuracy of the mechanical responses and the elastic properties of materials derived from convolutional neural networks or neural operators ~\cite{you2022, Mohammadzadeh2022, Kobeiss2022, Hoq2023}. However, the drawback of data-driven approaches is that they typically require a large number of datasets, on the order of thousands to millions of image pairs (\textit{i.e.,} an input image representing data to be analyzed and an output image corresponding to the expected outcome), in order for the machine learning model to learn their latent relation accurately. To this end, physics-informed machine learning~\cite{raissi2019,karniadakis2021physics} has emerged as an excellent alternative to traditional data-driven machine learning approaches to reduce the high demands on training data while maintaining high prediction accuracy~\cite{zhu2023cmame, yu2022cmame, lu2021_2}.

 Physics-informed neural networks (PINNs) have been utilized to determine the material parameters of homogenous structures within linear elastic, hyperelastic, and nonlinear-elastoplastic constitutive models~\cite{Nguyen-Thanh2020, Samaniego2020, haghighat2021, Wu2023}. Following the initial success of PINNs, this approach has further been extended to identify the material parameters of structures with additional internal inclusions~\cite{zhang2022, henkes2022, Hamel2023}, as well as heterogeneous linear elastic materials~\cite{kamali2023, chen2023}. At times, PINNs may struggle to converge in the presence of high-frequency functions or highly heterogeneous features~\cite{Wang2021,hao2023pinnacle,rathore2024challenges}. As such, the application of Fourier-feature PINNs in the context of system-biology-informed ordinary differential equations (ODEs) has also been explored~\cite{yazdani2020systems,Daneker2023}. 

The present research aims to tackle a complex inverse problem that involves realistic tissue structural patterns and nonlinear material constitutive models. In particular, we aim to investigate the applicability and generality of PINNs for characterizing the nonlinear, hyperelastic properties of biological tissues such as brain and tricuspid valve tissues, with an emphasis on discovering a high-fidelity map of complex elasticity fields in soft tissues, as well as investigating the relationship between neural network architectures and prediction accuracy. This foundational work will be extensible to accurately estimate \textit{in vivo} tissue properties from medical images, where the full-field strain will be used as reference data, and study the evolution of tissue microscopic elastic properties in response to mechanical loadings. We summarize the major contributions and significance of the current work as follows.
\begin{itemize}
  \item We introduce a physics-informed machine learning approach to discovering the elastic modulus distribution in complex heterogeneous hyperelastic materials that have traditionally been challenging and mathematically ill-posed.  
  \item Our improved neural network architecture consistently produces accurate estimates of full-field elasticity maps of complex material across multiple material constitutive models and biaxial stretch conditions.
  \item We demonstrated that PINNs have remarkable potential for providing highly accurate parameter estimations using just one data sample, even when there is a high level of noise in training data.
  \item This work provides a promising approach to studying the effect of material heterogeneity in their macroscopic behaviors and the evolution of material micromechanical properties in response to mechanical load. This is highly applicable to studying tissue growth and remodeling~\cite{myers2014, Oomen2018, fruleux2019}, and detecting early mechanical markers of cancer cell invasion~\cite{Hormuth2019}.
\end{itemize} 

\section{Results}
\label{sec:results}
Understanding how tissue microstructure relates to its mechanical properties and behaviors is crucial for maintaining tissue health. While some studies have explored how tissue microstructure affects its macroscopic mechanical properties~\cite{Molladavoodi2013, Budday2020, Zhang2022jmbbm}, research on the relationship between tissue microstructure and the corresponding elasticity field is limited. This study aims to uncover the elasticity field within soft tissue microstructural components using PINNs. 

\subsection{Problem setup}
\label{re:setup}
\begin{figure}[!h]
\centering
\vspace{-2ex}
\includegraphics[width=1\textwidth]{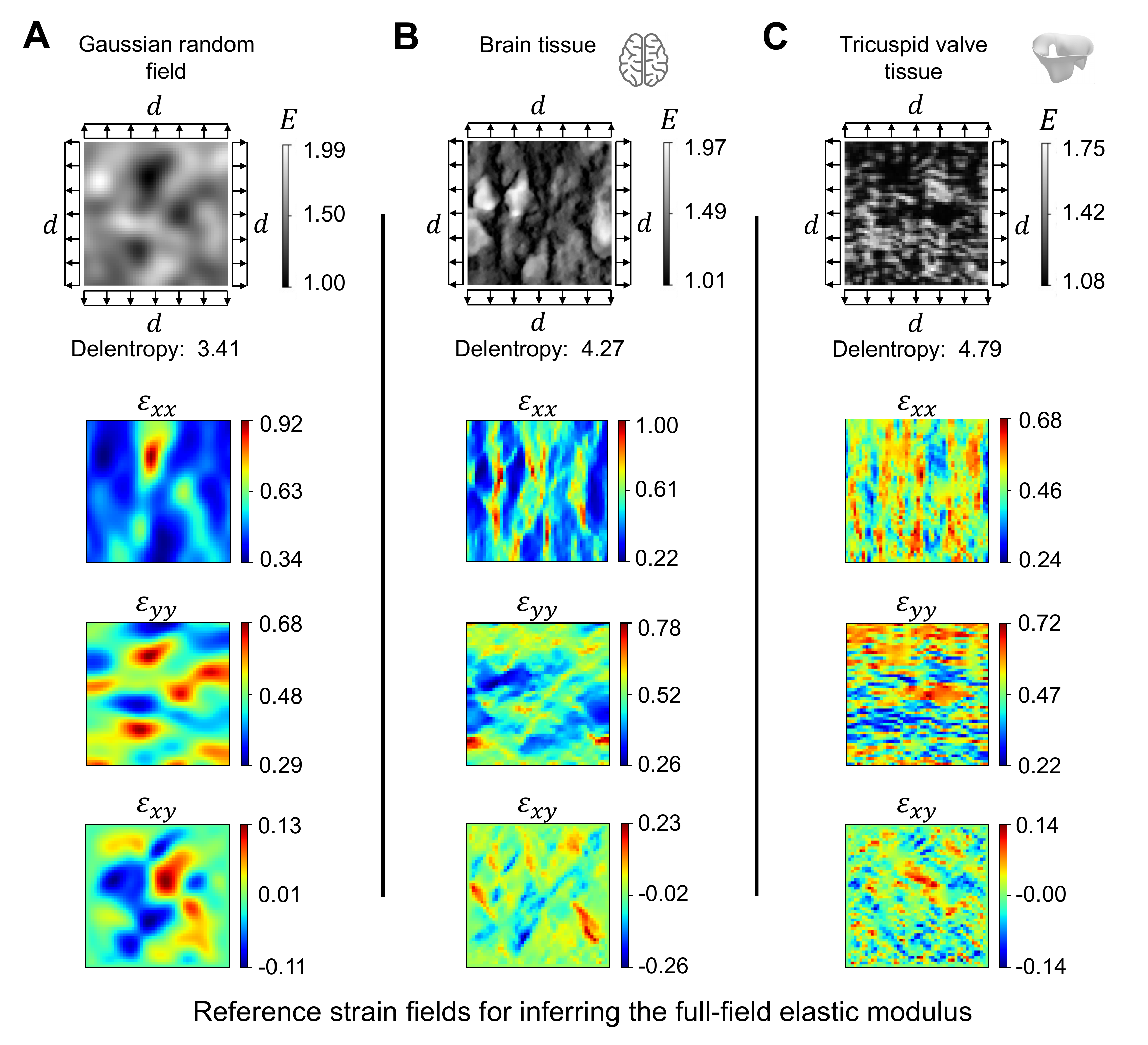}
\vspace{-5ex}
\caption{\textbf{Heterogeneous material examples.} The examples are presented in the undeformed configuration of unit square specimens. These examples represent the unknown elastic properties distribution of soft tissue that we seek to identify using PINNs. The structural complexity of the sample was quantified using Delentropy; higher values indicate greater complexity. (\textbf{A}) The structural pattern is synthetically generated using a GRF~\cite{Moore2015}. (\textbf{B}) The representation of the brain tissue microstructure was adapted from Koos et al.~\cite{Koos2016}. (\textbf{C}) The fibrous network of the tricuspid valve was adapted from Weinberg et al.~\cite{Weinberg2005}. Equibiaxial deformation of magnitude $d$ was applied to the undeformed configuration using FEA to obtain the ground truth strain distribution. We prescribed 40\% equibiaxial stretch to the specimen (\textit{i.e.,} 20\% to each side of the boundary) and estimated the resulting Green-Lagrangian strain using FEniCS~\cite{Anders2010}. The strain data were used as ground truth reference data to identify the ``unknown'' elastic moduli in PINNs inversely.} \label{tissue_examples}
\vspace{-5ex}
\end{figure}

Fig.~\ref{tissue_examples} illustrates the structural pattern of three heterogeneous tissue specimens considered in the present work. The tissue specimens are subjected to equibiaxial stretch as it is one of the standard experimental techniques for studying the structural-mechanical relationship of soft tissues, including the brain, heart valve, and skin tissues~\cite{Sacks2001, Labus2016, Meador2022}. The structural pattern of the first example was constructed using a Gaussian random field (GRF)~\cite{Moore2015} (Fig.~\ref{tissue_examples}A). The second example represents the cell bodies within a brain tissue~\cite{Koos2016} (Fig.~\ref{tissue_examples}B). The third example demonstrates the fibrous network of the tricuspid valve leaflets~\cite{Weinberg2005} (Fig.~\ref{tissue_examples}C). We adopted a complexity measure based on image information theory, namely Delentropy~\cite{larkin2016, Khan2022}, to evaluate the complexity of tissue microstructural patterns. Delentropy captures the relative relation between the local and global features by analyzing the gradient vector field of the image. This measure provides a useful metric to evaluate the performance of machine learning models in the context of data complexity.

To generate reference strain data, we first assumed an elastic moduli distribution ($E$) for the tissue microstructure in the image.  We then applied equibiaxial stretch to the external boundary of the model and computed the full-field strain data ($\epsilon_{xx}$, $\epsilon_{yy}$, and $\epsilon_{xy}$) using finite element analysis (FEA). Lastly, we integrated the full-field strain data in PINNs to guide the discovery of the heterogeneous elastic moduli distribution of the tissue microstructure, cell bodies, and fibrous network in PINNs. Detailed description on the calculation of the elasticity map is discussed in Section~\ref{data_generation}. In the first part of the study (Section~\ref{re:comparison}), we focused on identifying a suitable network architecture for identifying complex heterogeneous elasticity map. We applied 20\% of outward prescribed displacements to the external boundary of the square specimens and assumed a plane strain, compressible Neo-Hookean material model. In the second part of our study (Section~\ref{extended}), we expanded our explorations to evaluate the performance of the selected PINN architecture under various material models, levels of equibiaxial stretch, and noise levels in the strain data.

\subsection{PINN overview}
\label{re:pinn}
\begin{figure}[!h]
\centering
\vspace{-2ex}
\includegraphics[width=1\textwidth]{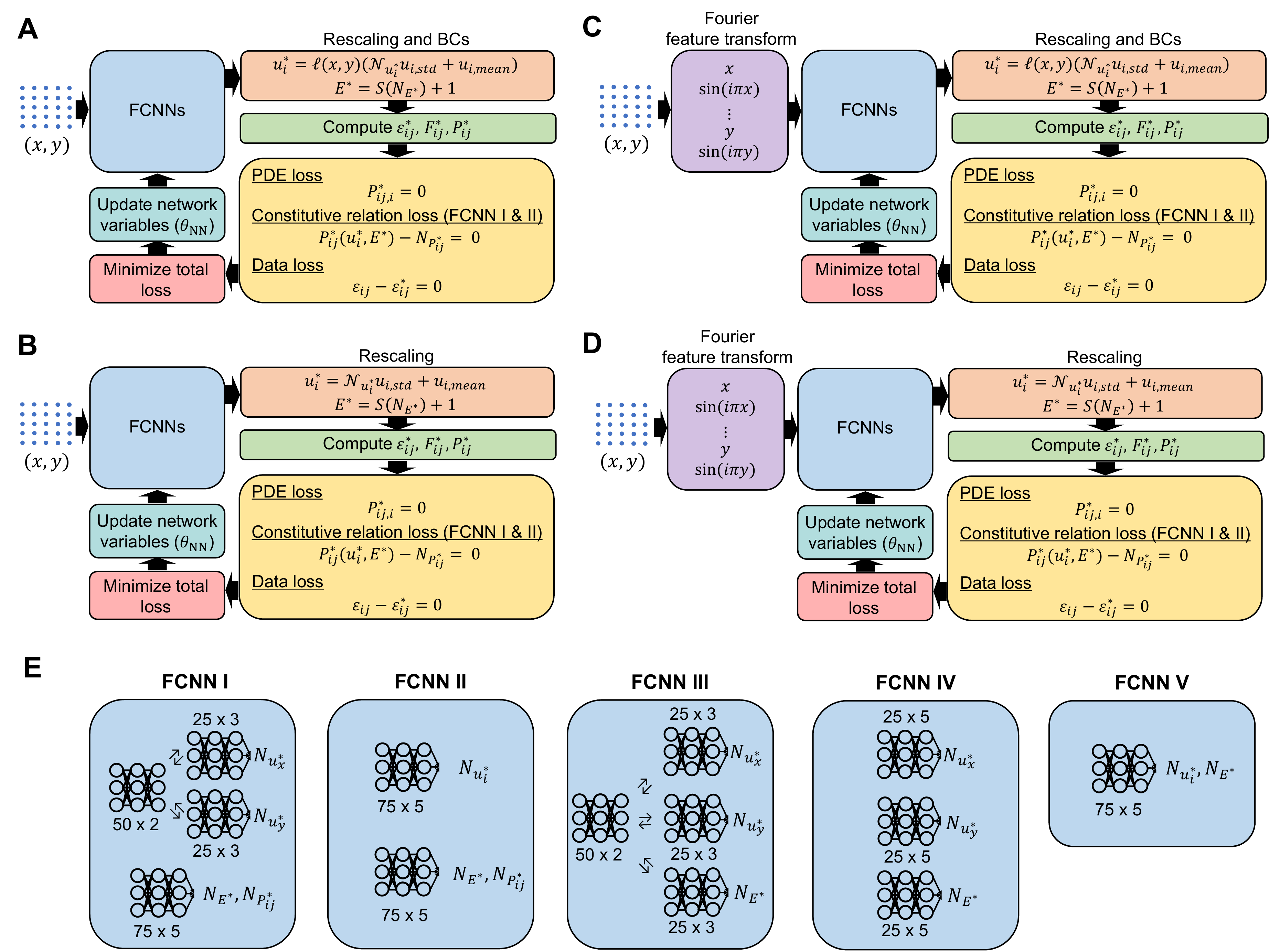}
\vspace{-2ex}
\caption{\textbf{Network architecture overview.} We established four PINN variants to examine the effectiveness of PINNs in identifying the heterogeneous elastic properties of hyperelastic materials. (\textbf{A}) A standard PINN with boundary conditions embedded as hard constraints. (\textbf{B}) A standard PINN without enforcing boundary conditions. (\textbf{C}) A Fourier-feature PINN with boundary conditions embedded as hard constraints. (\textbf{D}) A Fourier-feature PINN without enforcing boundary conditions. (\textbf{E}) For each PINN variant, we experimented with five types of FCNN architectures to determine the best arrangement of independent neural networks to achieve accurate estimates of the elastic properties distribution in complex materials. In Fourier-feature PINNs, the input coordinates of the PDE training points are mapped to a total of 10 Fourier modes (five Fourier modes for $x$ and five Fourier modes for $y$).}\label{network_arch}
\vspace{-2ex}
\end{figure}

In the current study, we evaluated the performance of four variations of PINNs, as shown in Figs.~\ref{network_arch}A to D. Figs.~\ref{network_arch}A and B represent standard PINNs architectures with and without boundary constraints, respectively. On the other hand, Figs.~\ref{network_arch}C and D display Fourier-feature PINNs architectures with and without boundary constraints, respectively. In addition, we examined the effect of five unique fully connected NN (FCNN) architectures (Fig.~\ref{network_arch}E) on the prediction accuracy of the constitutive elastic properties within each PINN variant. All simulations were performed on NVIDIA A100-SXM4-40GB GPUs. Detailed information related to Fourier-feature PINNs is described in Section~\ref{fourier_pinns}.

\subsubsection{PINN architectures}
In the following, we briefly summarize the five FCNN architectures (Fig.~\ref{network_arch}E). The output variables of FCNN I and FCNN II are $N_{u^*_{i}}$, $N_{P^*_{ij}}$, and $N_{E^*}$. Here, the output $N_{P^*_{ij}}$ is used to govern the constitutive relation loss between the stress tensor $P^*_{ij}(u^*_{i}, E^*)$, computed from the estimated $u^*_{i}$ and $E^*$, and the direct stress tensor, $N_{P^*_{ij}}$, output from the network. On the other hand, FCNN III, IV, and V are designed only to estimate $N_{u^*_{i}}$ and $N_{E^*}$ without ensuring a balanced constitutive relation loss.

To provide more details into FCNN I and II, FCNN I comprises two subnetworks. The first subnetwork has a shared network of two hidden layers with 50 neurons per layer. The shared network is split into two independent networks, each composed of two hidden layers with 50 neurons per layer. Each of the subnetworks predicts one of the displacement components (\textit{i.e.,} $N_{u^*_{x}}$ and $N_{u^*_{y}}$). The second subnetwork is a shared network consisting of five hidden layers with 75 neurons per layer that outputs $N_{E^*}$ and $N_{P^*_{ij}}$. In FCNN II, we similarly formulated two independent networks, each consisting of five hidden layers with 75 neurons per layer. The first independent network outputs the displacement components. In contrast, the second outputs the elastic moduli and the stress components, as the elastic moduli and stresses would share similar features given their direct mathematical relationship. 

Regarding FCNN III, IV, and V, in FCNN III, a two-hidden layer shared network is established to connect the displacement components and elastic moduli. The shared network is subsequently split into three independent networks with three hidden layers with 25 neurons per layer to individual features in the displacement component and elastic moduli. FCNN IV consists of three independent networks, each with five hidden layers and 25 neurons per layer. FCNN V represents the simplest form of neural network architecture, in which one shared network is formulated to estimate $N_{u^*_{i}}$ and $N_{E^*}$. 

The strong form of linear momentum equations, $P_{ij, i} = 0$, is considered as the governing equation of the physical system, wherein $P$ is the first Piola-Kirchhoff stress described in Section.~\ref{method:mat_models}. 

\subsubsection{Fourier-feature PINNs}
\label{fourier_pinns}
Standard PINNs often experience difficulties learning the parameter distribution of the materials with complex and random structural patterns. Seminal works have demonstrated the effectiveness of Fourier-feature transforms in solving benchmark PDEs problems and complex ordinary differential equations~\cite{Wang2021, Daneker2023}. As such, we extended our standard PINN architecture with Fourier-feature transforms to decompose the complex structural features into multiple layers to learn the high-frequency structural features, aiming to capture the intricate patterns through a multi-scale approach. 

We transformed the input coordinates through a Fourier-feature mapping with five Fourier modes per input variable before feeding it to the fully connected feedforward network, as depicted in Fig.~\ref{network_arch}. In particular, we retained the spatial coordinates, $x$ and $y$, with additional terms of $sin(i\pi x)$ and  $sin(i\pi y)$ for $i = 1, 2, \dots, 5$. Incorporating these features into the network can facilitate the learning of similar features and thus may reduce the computational difficulties for feature extraction. These very general Fourier features were found to provide good results. Additional hyperparameter fine-tuning could be applied for better results, such as $sin(\pi xy)$, but was not utilized in these cases for simplicity and proof of concept.

\subsection{Prediction accuracies of PINN architectures on elastic parameters estimation}
\label{re:comparison}
\begin{figure}[!h]
\centering
\vspace{-2ex}
\includegraphics[width=1\textwidth]{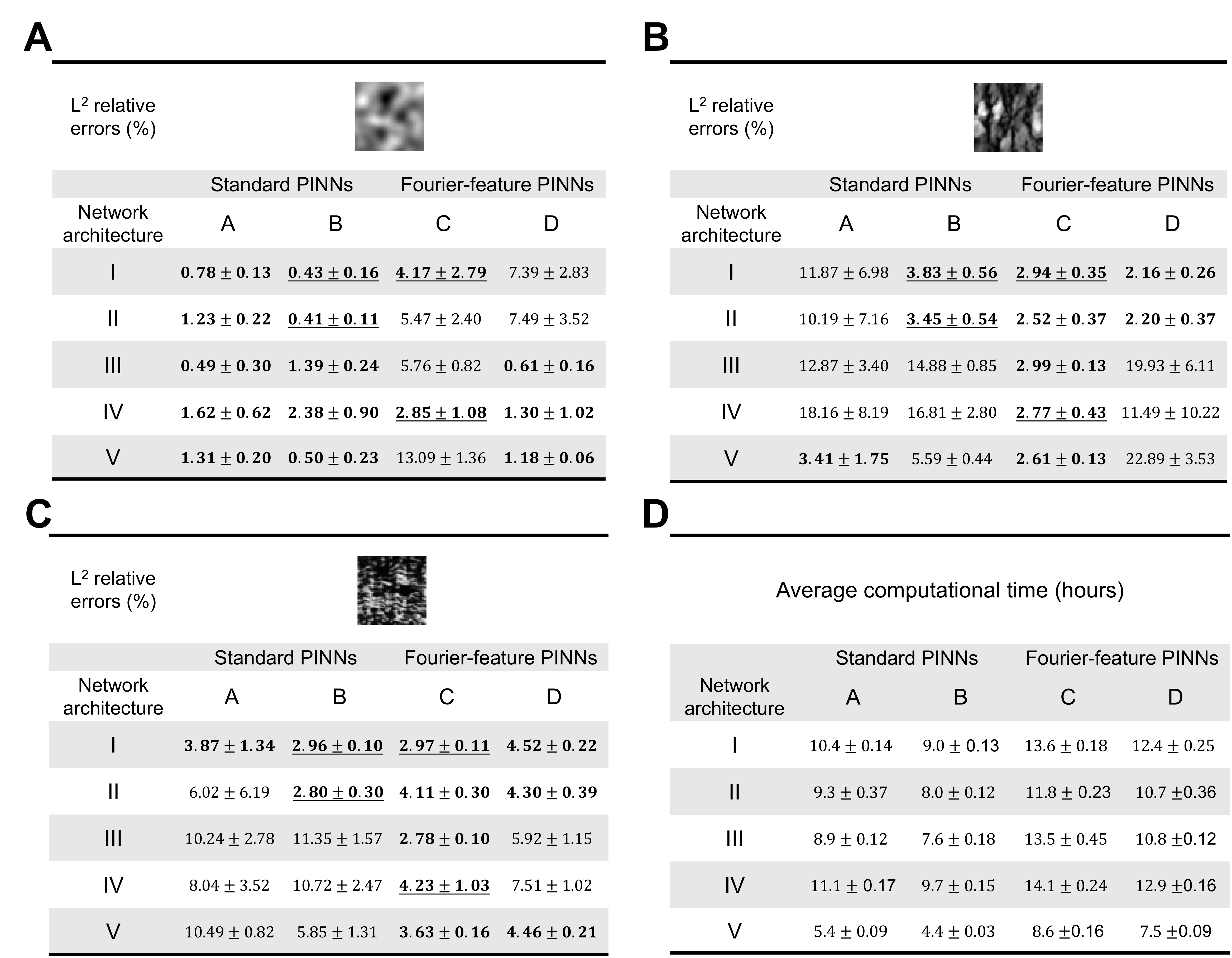}
\vspace{-2ex}
\caption{\textbf{Effect of network architecture on PINN prediction.} We conducted five evaluations for each PINN and FCNN architecture variant. (\textbf{A}--\textbf{C}) The mean and standard deviations of $L^2$ relative errors of the PINN-estimated full-field elastic moduli. The network architectures that produced less than 5\% $L^2$ relative errors are highlighted in bold. The architectures that consistently met the accuracy requirements across all three examples are underlined. Among the three network architectures that demonstrated high estimation accuracy across the examples, network architecture IIB was the most computationally efficient. (\textbf{D}) The computational time for each PINN architecture.}\label{overall_results}
\vspace{-2ex}
\end{figure}

This section discusses the prediction accuracies of the proposed architectures in estimating the complex elastic property distributions of soft materials that undergo large deformation. The material constitutive model was assumed plane strain compressible Neo-Hookean with Poisson's ratio $\nu=0.3$. The means and standard deviations of the $L^2$ relative errors of the estimated elastic modulus distribution of each example are shown in Figs.~\ref{overall_results}A to C. Each analysis was repeated five times with varying random seeds to ensure solution consistency. The FCNN that produced less than 5\% average $L^2$ relative errors are highlighted in bold. Furthermore, the architectures that consistently met the accuracy requirements across all three examples are underlined. The average computation time for each network architecture across the three examples and the five repeated tests are reported in Fig.~\ref{overall_results}D. 

In the first example (Fig.~\ref{overall_results}A), with standard PINNs, models A and B produced similar prediction accuracy across all five FCNN architectures. The difference in $L^2$ relative error between models A and B was less than 1\% across all FCNN architectures. On the contrary, with Fourier-feature PINN, FCNN III, IV, and V had higher $L^2$ relative errors in model C than D. These results suggested that, in the GRF example, Fourier-feature PINNs provided more accurate estimates of the elastic parameter distribution when the boundary conditions were unconstrained and constitutive law was unenforced. 

Similar $L^2$ relative error patterns were observed in the second and third examples (Figs.~\ref{overall_results}B and ~\ref{overall_results}C). We found that Fourier-feature PINNs outperformed standard PINNs in these two examples. Specifically, Fourier-feature model C yielded highly accurate estimates of the elastic parameter distribution with $L^2$ errors less than 5\% across all FCNN architectures. Neural network architectures that enforce constitutive relations (FCNN I and II) demonstrated superior prediction accuracy in both standard PINNs and Fourier-feature PINNs when compared to those that do not enforce constitutive relations (FCNN III, IV, and V). This observation highlights the importance of considering the underlying physics and problem characteristics when designing and training PINNs.

While we expected the network prediction accuracy would improve with additional information on the boundary. However, our experiments showed that embedding boundary conditions in the inverse analysis did not significantly improve the prediction of the unknown material parameters across the three examples in the standard PINNs. Conversely, enforcing boundary constraints substantially improved the estimation accuracy in the second example using Fourier-feature PINNs while producing comparable prediction accuracy in the first and third examples. These unexpected results could be because the network has learned the deformation boundary conditions from the provided strain reference data. Therefore, enforcing the boundary may not provide additional information to guide the optimization.

It is worth noting that Fourier-feature PINNs resulted in better estimates of the elasticity map in the second and third examples, but poorer estimates in the first example. This contradicting behavior suggests the existence of inherent trigonometric features in the full-field elasticity maps in the second and third examples, which were visibly apparent in the second and more subtle in the third. Therefore, the Fourier-feature transform provided additional guidance to the training process, subsequently improving the prediction accuracy. However, the first example did not seem to contain the same features. As a result, the Fourier-feature transform had a negative impact on the estimation of elastic property distribution.

Out of the 20 network architectures tested, the network architectures IB, IIB, IC, and IVC were found to produce highly accurate estimations consistently across the three examples. The results suggest that partitioning output variables into separate independent networks, based on their shared features, may substantially enhance prediction accuracy. Among these selected models, network architecture IIB provided accurate estimates of the elastic moduli, with average $L^2$ relative errors of 0.41\%, 3.45\%, and 2.80\% for the three examples respectively, while also being the fastest in terms of computational time.

In terms of computational time, we observed that Fourier-PINNs require a longer duration of training in comparison to their standard PINN counterparts. This can be attributed to the introduction of additional neural network parameters through the Fourier-feature mapping. Additionally, enforcing hard boundary conditions led to longer computational time in both standard and Fourier-feature PINNs.

\subsubsection{Comparison with the traditional inverse FEA method}

Here, we provide a comparative study between PINNs and inverse FEA for the GRF example. In inverse FEA, we integrated a topology optimization method, the adjoint method~\cite{Errico1997}, with FEA to minimize the $L^2$ relative errors between the ground truth and estimated strain tensor, subject to the constraint $1 \leq E \leq 5$, the same as that in PINNs. Within this approach, the PDE and boundary condition losses were naturally optimized within each FEA iteration. The adjoint inverse analysis was performed using the Dolfin-adjoint software~\cite{Mitusch2019} with an Interior Point Optimizer~\cite{Wachter2006}. The initial elasticity map was randomly generated. The maximum optimization iteration was set to 100, and the acceptable tolerance of the $L^2$ relative error for strains was set to $10^{-3}$. The results of our experiments are detailed in Appendix A. We found that the estimation accuracy using this optimization method is highly sensitive to the initial elasticity map.

\begin{figure}[!h]
\centering
\vspace{-2ex}
\includegraphics[width=1\textwidth]{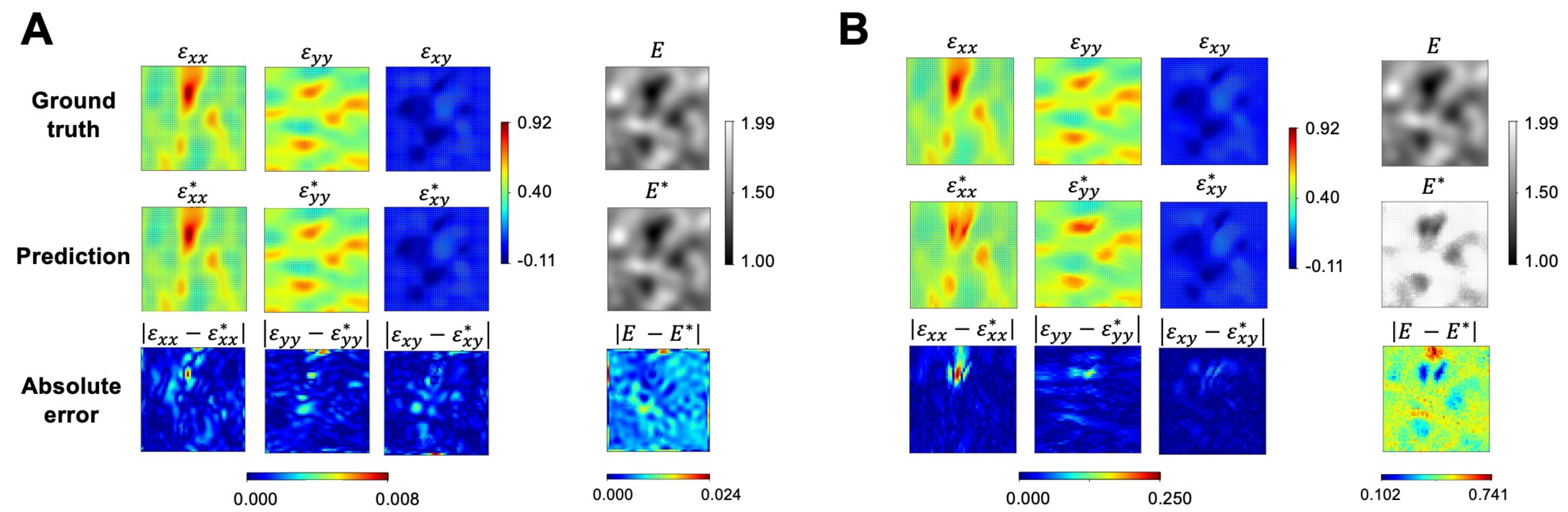}
\vspace{-2ex}
\caption{\textbf{Comparison of estimation accuracy using PINNs and the adjoint method.} (\textbf{A}) In the results from PINNs, the maximum absolute point-wise errors of strains and the full-field elasticity map were 0.008 and 0.024. The $L^2$ relative error of the elasticity maps was 0.56\%. (\textbf{B}) In the results from the adjoint method, the maximum absolute point-wise errors of strains and the full-field elasticity map were 0.25 and 0.741. The $L^2$ relative error of the elasticity maps was 29.22\%. Our proposed PINN architecture provides significant improvement in predicting the GRF heterogeneous field compared to the adjoint method.}\label{method_comp}
\vspace{-2.5ex}
\end{figure}

The estimated strain tensor and full-field elasticity map using PINNs with network architecture IIB and the best elasticity map estimation from the adjoint method are reported in Fig.~\ref{method_comp}. Fig.~\ref{method_comp}A presents results from PINNs, while Fig.~\ref{method_comp}B demonstrates results using inverse FEA. As shown, the maximum absolute point-wise error of strains resulting from PINNs was two orders of magnitude lower than that from the adjoint method, and the full-field elasticity map was one order of magnitude lower. The $L^2$ relative error of the full-field elasticity map estimated by PINNs was 0.56\%, while that estimated by the adjoint method was 29.22\%. Note that, in the adjoint method, the optimization process terminated early due to solution divergence in the FEA.

\newpage
\subsubsection{Convergence behavior of PINNs}
\label{re:convergences}
\begin{figure}[!h]
\centering
\vspace{-2ex}
\includegraphics[width=1\textwidth]{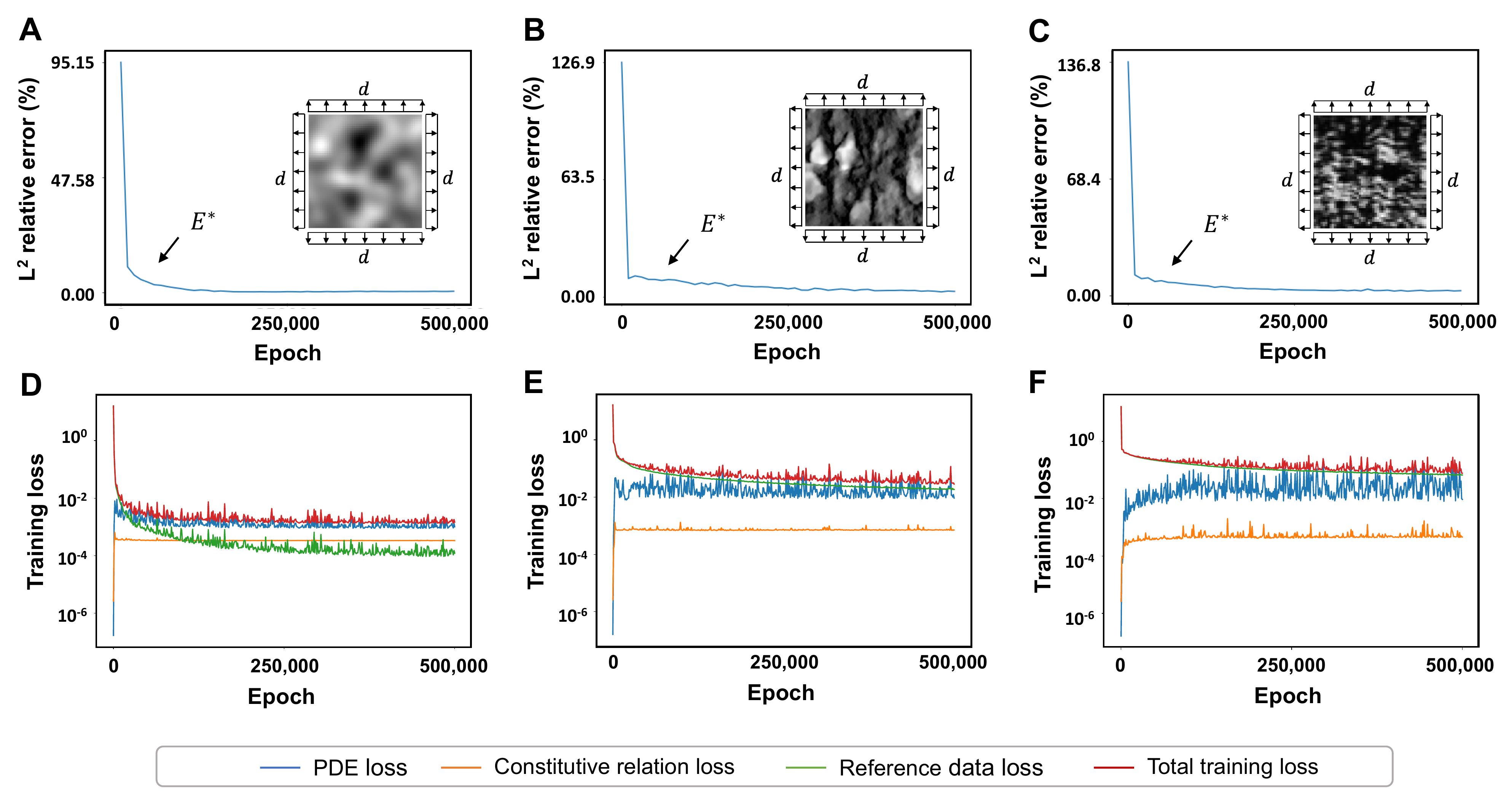}
\vspace{-2ex}
\caption{\textbf{PINN training behavior.} The $L^2$ relative errors of the full-field elastic parameters and the training loss convergence history are presented. We observed rapid convergences in the $L^2$ relative errors across all three examples. The $L^2$ relative errors of the elastic parameter predictions reduced below 5\% within 4,000 epochs in the first example, 250,000 epochs in the second example, and 150,000 epochs in the third example. The constitutive relation and PDE losses reach plateaus almost immediately at the beginning of the training across all three examples, while the total training loss follows the same decay rate as the reference data loss.}\label{cvg}
\vspace{-2.5ex}
\end{figure}

Fig.~\ref{cvg} shows the convergence histories of the elastic moduli and the training loss of each example using network architecture IIB. As shown, the training procedures and the estimated elastic moduli converge rapidly within the user-defined duration of $5\times10^5$ epochs. Figs.~\ref{cvg}A and D demonstrate the convergence evolution of the training procedure for the first example. The $L^2$ relative error of the elastic moduli first decays to below 5\% at around 4,000 iterations. The estimated elastic moduli continue to approach the FEA ground truth for the remaining iterations, reaching a final relative error of 0.56\%. Similarly, the total training loss decreases from $\mathcal{O}(10^{-1})$ to $\mathcal{O}(10^{-3})$ within the first 150,000 iterations. The PDE and constitutive relation losses have reached a plateau immediately at around 1,000 iterations and consistently maintain $\mathcal{O}(10^{-3})$ and $\mathcal{O}(10^{-4})$ relative error. Similar patterns are observed in the second example (Figs.~\ref{cvg}B and E) and C (Figs.~\ref{cvg}C and F). The $L^2$ relative error of the elastic moduli in the second example decays to below 5\% at around 250,000 iterations and steadily reduces the error to 2.91\% at the end of the training. The $L^2$ relative error of the elastic moduli in example C decays to below 5\% at around 150,000 iterations and reduces the error to 2.78\% at the end of the training. Across the three examples, total training loss converging behavior aligns more with data loss, suggesting that the training procedure is more sensitive to deviation in the strain reference data. The estimations and pointwise errors of displacements, strains, and elasticity maps resulting from network IIB are detailed in Appendix B.

\newpage
\subsection{Extended experiments on selected PINN architectures}
\label{extended}
In Section~\ref{re:comparison}, we identified four accurate network architectures that consistently produced accurate estimates of the elasticity map across the test examples. In this section, we first identify the best-performing network architecture by evaluating its estimation accuracy on various heterogeneous tissue patterns. Afterward, we examine the robustness of the selected PINN architecture in the presence of noisy strain data, as well as when extrapolated to different mechanical loading conditions and material constitutive models.

\subsubsection{Hetergeneous tissue patterns from ultrasound images}

We identified that networks IB, IIB, IC, and IVC consistently yielded accurate full-field elastic parameter prediction of the tissue with less than 5\% $L^2$ relative error in Section~\ref{re:comparison}. Among them, IB and IIB were standard PINNs, while IC and IVC were Fourier-feature PINNs. In the following, we apply these four network architectures to analyze three cases of malignant breast cancer tissues captured from ultrasound images from a public dataset~\cite{Al-Dhabyani2020}. It is important to note that we used the same hyperparameters as in the previous examples to demonstrate consistency in our methods. However, the hyperparameters could be further refined to achieve better estimation results.

In the three breast cancer tissue examples, we generated the strain data following the same forward FEA procedures as the brain and tricuspid valve tissue examples. In the forward FEA, 40\% equibiaxial stretch was applied to the boundary of the specimens. Similarly, the Neohookean hyperelastic material model with Poisson's ratio of 0.3 was assumed. The full-field elastic parameters were estimated based on the pixel values of the images as described in Section~\ref{data_generation}.

\begin{figure}[!h]
\centering
\vspace{-2ex}
\includegraphics[width=1\textwidth]{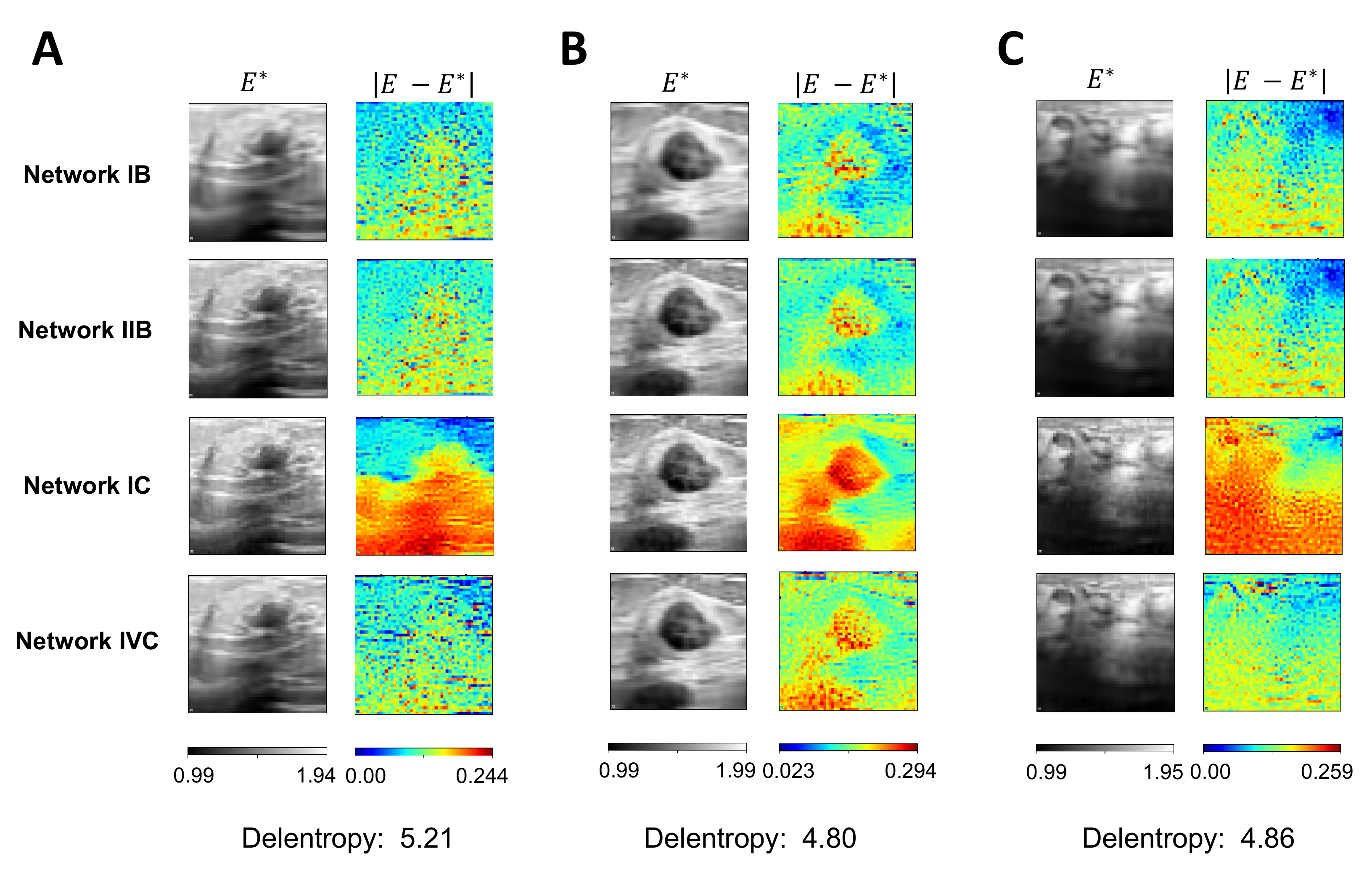}
\vspace{-2ex}
\caption{\textbf{Breast cancer tissues.} Three examples of malignant breast cancer tissues were chosen from a public dataset~\cite{Al-Dhabyani2020} to demonstrate the applicability of our proposed PINN architectures to a wider range of tissue patterns with delentropy complexity values beyond those in the previous examples. The pointwise absolute errors were all under 0.3 across the three tissue examples and network types. In example (\textbf{A}), the $L^2$ relative errors for networks IB, IIB, IC, and IVC were found to be 5.60\%, 5.20\%, 10.24\%, and 3.13\%, respectively. Meanwhile, the errors were 12.01\%, 11.51\%, 13.05\%, and 8.48\%, respectively in example (\textbf{B}) and 6.92\%, 6.69\%, 11.80\%, and 4.57\%, respectively in example (\textbf{C}). Network IVC consistently achieved less than 10\% $L^2$ relative errors in estimating heterogeneous elasticity maps.}\label{cancerTissues} 
\vspace{-2.5ex}
\end{figure}

The full-field elastic parameter estimations and pointwise errors of the breast cancer tissues are shown in Fig.~\ref{cancerTissues}. These examples have delentropy complexity values ranging from 4.80 to 5.21, which are larger than the values of the previous examples. In cancer tissue example A, the $L^2$ relative errors for networks IB, IIB, IC, and IVC were 5.60\%, 5.20\%, 10.24\%, and 3.13\%, respectively. For example B, the errors were 12.01\%, 11.51\%, 13.05\%, and 8.48\%, respectively. In example C, the errors were 6.92\%, 6.69\%, 11.80\%, and 4.57\%, respectively. Across the three examples, network IVC consistently produced accurate prediction with less than 10\% $L^2$ relative errors, while network IIB consistently delivered the best estimations of full-field elasticity maps within standard PINN architectures. Due to the lower computational time required, network IIB was used in all remaining analyses in the present work.

\subsubsection{Effect of noise in reference data}
\begin{figure}[!h]
\centering
\vspace{-2ex}
\includegraphics[width=1\textwidth]{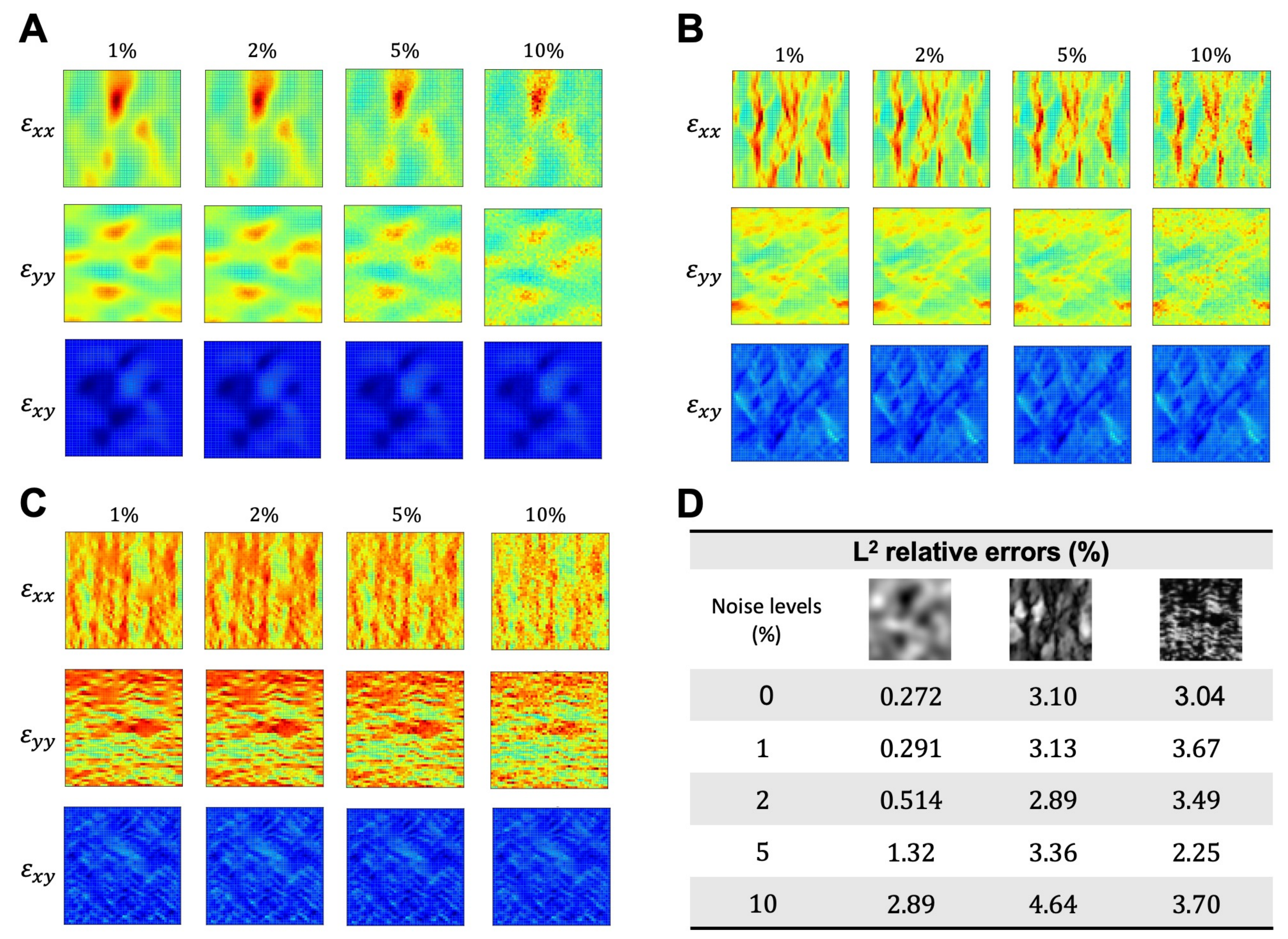} 
\vspace{-2ex}
\caption{\textbf{Effect of noisy reference data.} (\textbf{A-C}) The reference strain data of the GRF, brain tissue, and tricuspid valve tissues subjected to $1\%$, $2\%$, $5\%$, and $10\%$ white Gaussian noise. (\textbf{D}) The $L^2$ relative errors are below $5\%$ across all examples and noise levels.}\label{noise_results} 
\vspace{-2ex}
\end{figure}

Strain data, either from experimental measurements or strain elastography, are typically subject to varying degrees of noise that could interfere with the accuracy of inverse estimation. To assess the robustness of the proposed network architecture, we injected $1\%$, $2\%$, $5\%$, and $10\%$ of white Gaussian noise into the strain tensor of each example. The reference noisy strain data, along with the $L^2$ relative errors of the estimated elasticity maps are shown in Fig.~\ref{noise_results}. Our experimentation revealed that the estimation accuracy of the elasticity map decreased as the level of noise increased in the first example. However, the $L^2$ relative errors were similar across different noise levels in the second and third examples. This impartial behavior towards noise measurements could be attributed to the inherently complex structural patterns within the brain and tricuspid valve tissues. Nonetheless, the selected network architecture was able to estimate the elasticity maps with a high degree of accuracy even when the reference data was corrupted with a high level of noise, as the $L^2$ relative errors consistently remain below $5\%$ across all examples and noise levels (Fig.~\ref{noise_results}D).

\subsubsection{Estimating heterogeneous material parameters of plane stress Neo-Hookean model subjected to varying levels of biaxial stretch}
\label{varying_stretch_study}
In the following experiments, we used the same full-field elasticity maps in Fig.~\ref{tissue_examples} for each example. We modified the material stress-strain behavior to a plane stress Neo-Hookean model and assumed a Poisson's ratio of 0.45. Each model was subjected to three levels of biaxial stretching ($10\%$, $20\%$, $40\%$) independently, and the reference strain data were obtained using FE simulation. Overall, Network architecture IIB accurately estimated full-field elasticity maps with $L^2$ relative errors within $5.04\%$ in all experiments.

\begin{figure}[htbp]
\centering
\vspace{-2ex}
\includegraphics[width=1\textwidth]{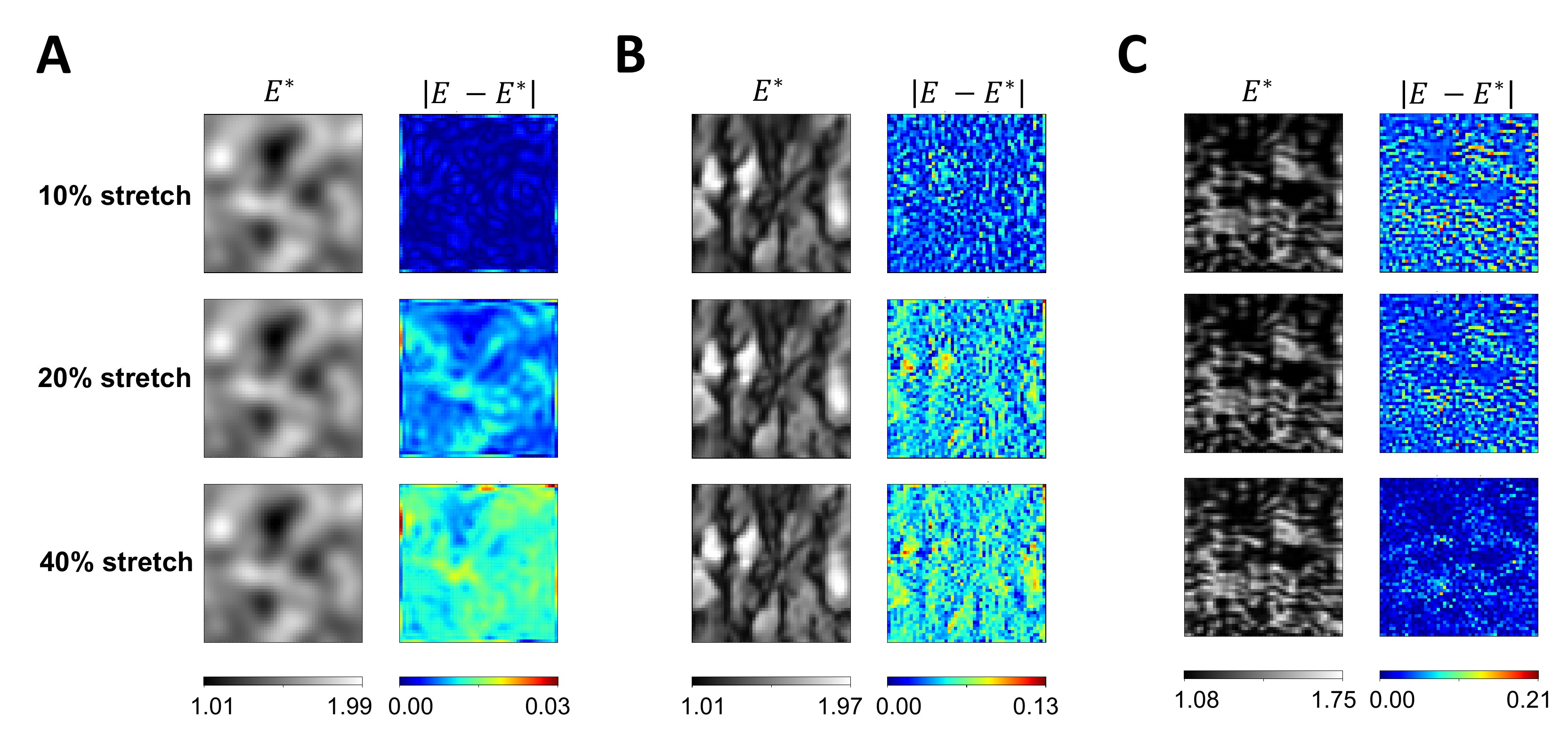}
\vspace{-2ex}
\caption{\textbf{Effect of biaxial stretch intensity on PINN predictions and pointwise errors.} (\textbf{A} and \textbf{B}) The pointwise error increased as the example underwent larger deformation due to higher prescribed biaxial stretch. On the contrary, (\textbf{C}), higher biaxial stretch produced a more accurate estimate of the full-field elasticity map. In general, the maximum pointwise error of the estimated full-field elasticity maps increased with a higher Delentropy value and complexity in the structural pattern of the examples.}\label{vary_stretch}
\vspace{-2ex}
\end{figure}

Fig.~\ref{vary_stretch} shows the estimated elasticity maps and pointwise errors of the three examples subject to varying levels of biaxial stretch. In the GRF example (Fig.~\ref{vary_stretch}A), it was observed that as the percentage of stretch increased, the $L^2$ relative errors for the elasticity map estimation also increased. The $L^2$ relative errors were found to be $0.17\%$, $0.63\%$, and $0.95\%$ for stretch levels of $10\%$, $20\%$, and $40\%$, respectively. Similar trends of the $L^2$ relative errors were observed in the brain tissue example presented in Fig.~\ref{vary_stretch}B. The $L^2$ relative errors were found to be $2.38\%$, $3.63\%$, and $4.01\%$ for stretch levels of $10\%$, $20\%$, and $40\%$, respectively. However, in the tricuspid valve tissue example (Fig.~\ref{vary_stretch}C), higher $L^2$ relative errors were found with a lower percentage of biaxial stretch. The $L^2$ relative errors were found to be $5.04\%$, $4.25\%$, and $2.23\%$ for stretch levels of $10\%$, $20\%$, and $40\%$, respectively. Moreover, the maximum pointwise error across the estimated full-field elasticity maps generally increased with a higher Delentropy value in the underlying tissue structures. The maximum pointwise error of the full-field elasticity maps across the three examples were 0.03, 0.13, and 0.21, respectively.

\subsubsection{Estimating heterogeneous material parameters of compressible Neo-Hookean, incompressible Mooney Rivlin, and incompressible Gent models}
In the following experiments, we modified the material constitutive model to plane stress Neo-Hookean~\cite{bonet_wood_2008}, incompressible Mooney Rivlin~\cite{Mooney1940}, and incompressible Gent~\cite{Gent1996}. Detailed formulations of these material models are described in Fig.~\ref{method:mat_models}. The full-field elasticity maps in Fig.~\ref{tissue_examples} were used to represent the unknown distribution of the material variable $E$, $\mu_1$, and $\mu$ in the Neo-Hookean, Mooney Rivlin, and Gent models, respectively. All examples were subjected to $40\%$ biaxial stretch, and the reference strain data were obtained using FE simulation. Network architecture IIB accurately estimated full-field material parameter maps with $L^2$ relative errors well within $5\%$ in all experiments, except the tricuspid valve example with Gent model had a $L^2$ relative error of $7.05\%$.

\begin{figure}[htbp]
\centering
\vspace{-2ex}
\includegraphics[width=1\textwidth]{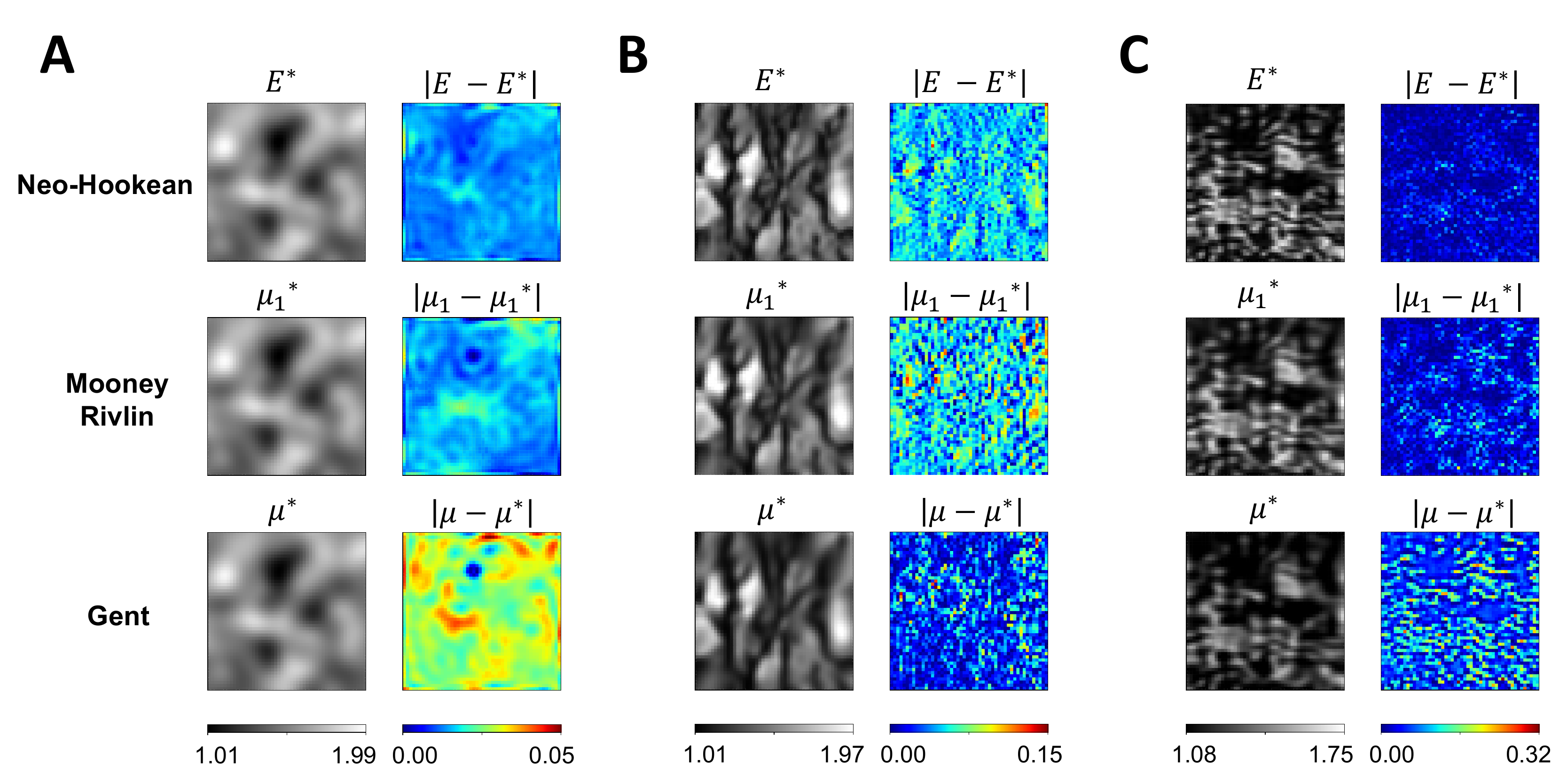}
\vspace{-2ex}
\caption{\textbf{Effect of material constitutive models on PINN predictions and pointwise errors.} (\textbf{A} and \textbf{C}) An increase in the complexity of the mathematical description of the material models correspondingly resulted in an increase in pointwise errors. In both examples, the Gent model reported the highest pointwise errors, while the Neo-Hookean model returned the lowest.  Conversely, in (\textbf{B}), the Gent Model had the lowest pointwise errors. Similar to Fig.~\ref{vary_stretch}, examples with higher structural complexity produced higher pointwise errors.}\label{vary_model}
\vspace{-2ex}
\end{figure}

Fig.~\ref{vary_model} shows the estimated parameter distribution and pointwise errors of the three examples characterized by plane stress compressible Neo-Hookean, incompressible Mooney Rivlin, and incompressible Gent models. In the GRF (Fig.~\ref{vary_model}A) and tricuspid valve (Fig.~\ref{vary_model}C) examples, it was observed that the $L^2$ relative errors of the heterogeneous material parameter map increased as the complexity of the material constitutive model increased. For the GRF example, the $L^2$ relative errors were found to be $0.95\%$, $1.08\%$, and $2.00\%$ for Neo-Hookean, Mooney Rivlin, and Gent, respectively. Similarly, for the tricuspid valve example, the $L^2$ relative errors were found to be $2.23\%$, $3.64\%$, and $7.05\%$ for Neo-Hookean, Mooney Rivlin, and Gent, respectively. On the contrary, in the brain tissue example (Fig.~\ref{vary_model}B), the estimated full-field material parameters of the Gent model were found to have the lowest $L^2$ relative error, which were observed to be $4.01\%$, $4.39\%$, and $2.70\%$ for Neo-Hookean, Mooney Rivlin, and Gent, respectively. Consistent with Section~\ref{varying_stretch_study}, higher Delentropy value in the complexity of the examples yielded a higher pointwise error of the estimated full-field material parameters. The maximum pointwise error of the full-field elasticity maps across the three examples were 0.05, 0.15, and 0.32, respectively.

\section{Discussion}
\label{sec:discussion}
In the current work, we examined the effects of PINNs for uncovering the heterogeneous material property of nonlinear elastic materials with complex tissue microstructures. We constructed 20 PINN architectures to identify the full-field elastic modulus of soft tissues characterized by an isotropic incompressible Neo-Hookean material model. We used the same hyperparameters (\textit{e.g.,} loss weights, learning decay rate, and optimization algorithm) in each analysis to ensure generality and that the estimated results are independent of parameter tuning. In our computational experiments, we observed that network architecture IIB consistently yielded accurate estimations of the full-field elastic modulus across a wide range of tissue structural patterns and simultaneously took the least computational time. These promising results led to further evaluation of the robustness and generality of architecture IIB on examples subjected to various biaxial stretch intensities and more complex material models, such as the incompressible Mooney Rivlin and Gent models. Across all experimentations, the selected PINN architecture consistently delivered highly accurate estimates of the full-field material constants even in the presence of noisy data. This study presents a critical first step toward a robust PINN model for identifying full-field heterogeneous material constants in soft tissues.

Our proposed approach extends on the works of Kamali et al., and Chen et al. ~\cite{kamali2023, chen2023}, which focused on identifying the linear elastic properties of heterogeneous materials under small deformation and infinitesimal strains. Here, we applied large deformation to our examples and used three hyperelastic material models to better characterize the material behavior of soft tissues. In addition, our methodology differed from that of Kamali et al.~\cite{kamali2023} as we utilized full-field strain computations from displacement output to ensure local kinematics were accurately satisfied, instead of using strain fields as network input with boundary tractions to optimize network parameters. Furthermore, our approach eliminates the necessity for stress data, which may not always be readily available and measurable, making it more adaptable to clinical imaging when boundary conditions are not precisely known. While the Poisson's ratio of soft tissue may vary between 0.3 to 0.5~\cite{Islam2020}, measuring the compressibility of soft tissue non-invasively is non-trivial. In fact, most clinical strain elastography assumes a constant Poisson's ratio. To that end, our study focused on identifying the heterogeneous stiffness maps. It is important to note that the current strain data were generated from simulated biaxial numerical experiments, whereas the boundary and loading conditions experienced by tissues in \textit{in vivo} data could differ depending on the type of examination. Therefore, additional validation of the present framework using medical images will be undertaken.
 
Our proposed PINN architecture may lead to new non-invasive approaches for quantifying the elasticity map utilizing the full-field strain measurement obtained from magnetic resonance imaging, optical coherence tomography, or strain-based ultrasound elastography. Shear wave elastography is a common imaging technique for detecting cancerous tissue by measuring tissue elasticity maps through shear wave velocity. However, commercially available US elastography systems assume material linearity in soft tissues, despite the high degree of nonlinearity and heterogeneity present~\cite{Ryu2017}. Our work demonstrated that our architecture can accurately estimate the elastic properties distribution using only strain as the training data. PINNs enable simultaneous estimation of the full-field stress and elastic moduli; this characteristic circumvents the need for homogeneous material assumption and may provide a more accurate quantitative measure of the stiffness map, therefore leading to a more accurate diagnosis. 

Our method provides a compact, flexible, unified framework for material properties studies. Our network architecture offers immense flexibility to integrate additional empirical data and relevant physical principles and knowledge into the network. For example, the current network could be extended to analyze the material property distributions of anisotropic materials by including additional invariant terms in the material constitutive models. Further, experimental data could also be easily integrated into the network to constrain the optimization path of the network variables, which in turn may further improve the elastic property prediction accuracy. While the current study focuses on the heterogeneous properties of soft materials, our architecture is highly suited for other engineering applications, with an easily extensible framework to include data such as the temperature field and thermodynamics conservation laws to study the material microstructure progress during crack propagation, helping accelerate progress in material designs and optimization.

While optimizing the computational time was not a focus within the scope of this work, there are a few approaches one could make to reduce the computational cost. First, our studies adopted the FE nodal coordinates (a total of 5101 nodes) as the PDE training points and reference data collocation points. One may experiment with the effect of reducing the number of training points on estimation accuracy~\cite{lu2021,Wu2023,wu2023comprehensive}. Second, as noted in Section~\ref{re:convergences}, the training loss appears to converge rapidly within the first 250,000 iterations of the analysis; depending on the level of accuracy desired, one may also reduce the number of iterations needed. Third, we used the Adam optimizer in our training; one may experiment with the use of second-order optimizers~\cite{Liu1989,rathore2024challenges} to achieve faster estimation convergence while maintaining the same order of accuracy. Lastly, there has been active research on methods to reduce GPU memory consumption and wall time in physics-informed machine learning, such as mixed precision~\cite{hayford2024speeding} which can reduce the computational cost by up to 50\%. 

In addition to discovering the full-field heterogeneous material parameters, advancement in machine learning has inspired complementary efforts to develop data-driven approaches for identifying material constitutive models appropriate for soft tissues~\cite{Linka2023}. By discovering constitutive models, we can gain insights into the homogenized, macroscopic material behaviors, providing a mathematical description to best describe phenomenological observation. By discovering full-field material parameters, we can better understand the internal structure and elastic properties of tissues meso- or microscopically, bringing new insights into the behavior of internal tissue structure in response to external mechanical forces that lead to the mechanical responses that were observed macroscopically. Robustly coupling the multiscale mechanical behaviors of soft tissues remains an open challenge~\cite{Guo2022}. Developing a multiscale machine learning model that couples soft tissue behaviors across macro-, meso-, and microscale to better inform the mechanobiology of soft tissues would be an exciting long-term objective to explore in future studies. 

In this context, it is essential to incorporate realistic microstructures derived from high-resolution imaging modalities such as MRI, OCT, or strain elastography. These techniques provide detailed data that can significantly enhance the accuracy and physiological relevance of our models. Integrating this real data will allow us to move beyond randomly distributed material parameters and towards more precise, data-driven representations of tissue microstructure. This approach will not only improve the fidelity of our models but also provide deeper insights into the behavior of biological tissues under various loading conditions, thus better informing the mechanobiology of soft tissues.

Lastly, in our work, we have chosen three examples with complex morphology to encapsulate the wide variety of structural patterns in soft tissues. Here, the structural complexity of the tissue sample was quantified using Delentropy. We hope this measure will help guide the designs of the FCNN architecture tailored to individual tissue specimens. It should be noted that the depth and width of FCNNs may need to be adjusted to accommodate the increasing complexity of the tissue structure. Further, the network output and governing equations can be easily extended to accommodate three-dimensional analysis. However, in the case of 3D analysis, the accuracy of estimation could be compromised if strain measurements are limited to the surface of a specimen. Fortunately, strain values will be accessible throughout the region of interest in strain elastography.

\section{Conclusions} 
\label{sec:conclusion}
In conclusion, we identified an effective PINN architecture, from 20 architecture variants, that consistently demonstrated excellent estimation accuracy of elasticity maps on three unique heterogeneous materials (GRF, brain tissue, and tricuspid valve tissues). The selected PINN model demonstrated robustness and high estimation accuracy on three different hyperelastic constitutive models (Neo-Hookean, Mooney Rivlin, and Gent), various levels of biaxial stretch, and up to 10\% noise in the training data. Further, the selected architecture demonstrated excellent extrapolation to additional soft tissue structures that have Delentropy complexity values beyond the initial experiment examples. The results of this work provide a promising framework for identifying the heterogenous elastic properties in soft materials from images, paving the way to accelerate research progress in material designs, material optimization, and medical applications. 

\section{Methods}
\label{sec:methods}
In this work, we investigated the accuracy of standard PINNs and Fourier-feature PINNs in estimating the material properties of heterogeneous materials in 2D. Each PINN method consisted of two variants: one with boundary conditions embedded as hard constraints and one without. Additionally, we established five network architectures with a unique arrangement of independent neural networks. The current study consists of two parts. First, a total of 20 neural network architectures were applied to each example to determine the most effective architecture for predicting the elastic modulus distribution in complex materials. Second, once a top-performing architecture was identified, it was modified to include noise in the training data and adapted to three different hyperelastic material models, and biaxial stretch levels to assess their performance when extrapolating other mechanical scenarios. Details on the data generation and details pertaining to the neural network architectures are provided in Sections \ref{data_generation} to \ref{pinns_framework}.  

\subsection{Reference data generation}
\label{data_generation}
In our experience, the displacement fields often do not provide sufficient local data to uncover the fine features in heterogeneous materials. Therefore, we generated strains ($\epsilon_{ij}$) from FEA to regulate the data loss and train our neural networks. To set up finite element models, we first prescribed heterogeneous patterns and elastic properties to a finite element mesh. In all three examples, the finite element mesh consisted of a unit square with 10,000 equilateral triangular elements and 5,101 nodes.

In the first example, the heterogeneous pattern was generated by assigning random intensity values to the nodal coordinates using a mean-zero GRF defined as 
\begin{equation*}\label{grfs}
f(x) \sim \mathcal{GP}(0, k_l(x_1, x_2)), 
\end{equation*}
where $k_l(x_1, x_2)$ is the Guassian kernal and $l$ is the correlation length. We chose a radial basis kernel function with $l=0.1$ in this work. In the second and third examples, we utilized the tissue microstructure images obtained from Koos et al.~\cite{Koos2016} and Weinberg et al.~\cite{Weinberg2005} to map the complex microstructure patterns to a finite element mesh. First, we cropped a section of the image and transformed the dimension to a unit square. Then, we used min-max normalization to scale the image intensity values to a range from 0 to 1. Finally, we interpolated the image intensity using cubic interpolation to map the heterogeneous feature of the brain and tricuspid valve leaflet tissue onto the finite element mesh.

In these examples, Young's modulus of each element was calculated based on the average nodal intensity value, defined as
\begin{equation*}
E_{ele} = \sum_{i=1}^3 f(x_i) + 1,
\end{equation*}
where $f(x_i)$ is the pixel intensity associated with node $i$ on the triangle element; a higher intensity value indicates greater stiffness. The finite element mesh was deformed by applying an equibiaxial displacement ($d = 0.2$) to the boundaries, as shown in Fig.~\ref{tissue_examples}. 

\subsection{Hyperelastic constitutive models}
\label{method:mat_models}
We evaluated the effectiveness of our proposed PINN architecture in identifying the heterogeneous properties of soft tissues on three popular hyperelastic constitutive models: the compressible Neo-Hookean~\cite{bonet_wood_2008}, incompressible Mooney Rivlin ~\cite{Mooney1940}, and incompressible Gent~\cite{Gent1996} models.

\textbf{Neo-Hookean model.} The Neo-Hookean model~\cite{bonet_wood_2008} is the simplest nonlinear elastic material model characterized by the 2D strain energy density function, $\Psi$,
\begin{equation*} 
\Psi(I_1, J) = \frac{1}{2}\lambda[log(J)]^2-\mu log(J)+ \frac{1}{2}\mu (I_1-2),
\end{equation*}
where $I_1$ is the first principal invariants denoted as $I_1 = \text{trace}(F^T\cdot F)$, $F$ is the deformation gradient denoted as $F_{ij} = \delta_{ij} + u_{i, j}$, $u_i$ is the displacement vector, and $\lambda$ and $\mu$ are the Lam\'e's elasticity parameters. The first Piola-Kirchhoff stress of hyperelastic material models is defined as the derivative of the strain energy function with respect to the deformation gradient. In Neo-Hookean model,  
\begin{equation*}
P = \frac{\partial \Psi}{\partial F} = \mu F+[\lambda log(J)-\mu]F^{-T}.
\end{equation*}
The plane-strain formulation of the Neo-Hookean model is expressed as
\begin{gather*}
\lambda = \frac{E\nu}{(1+\nu)(1-2\nu)}, \quad \textrm{and} \quad \mu = \frac{E}{2(1+\nu)}. 
\end{gather*}
In the case of plane stress, the Lam\'e constant $\lambda$ becomes 
\begin{equation*}
\lambda = \frac{2\lambda\mu}{\lambda + 2\mu}.
\end{equation*}

\textbf{Mooney Rivlin model.} The Mooney Rivlin model~\cite{Mooney1940} is made of a linear combination of the first and second principal invariants, $I_1$ and $I_2$. The strain energy density function is expressed as
\begin{equation*} 
\Psi(I_1, I_2) = \frac{\mu_1}{2}(I_1-2)+ \frac{\mu_2}{2}(I_2-2),
\end{equation*}
where $\mu_1$ and $\mu_2$ are empirically determined material constants. The second principal invariant is defined as $I_2 = \frac{1}{2}[I_1^2 - \text{trace}(C\cdot C)]$, where $C$ is the right Cauchy-Green tensor defined as $C=F^{T}\cdot F$. Typically, $\mu_2$ is much smaller than $\mu_1$. Further, for small deformation problems, the sum of $\mu_1$ and $\mu_2$ leads to the shear modulus $\mu$. The first Piola-Kirchhoff  stress is 
\begin{equation*}
P = \mu_1 F+ \mu_2 (I_1F-F \cdot C).
\end{equation*}

\textbf{Gent model.} The Gent model~\cite{Gent1996} captures material elasticity using logarithm functions. The strain energy density function is 
\begin{equation*} 
\Psi(I_1) = -\frac{\mu J_m}{2}\ln(1-\frac{I_1-2}{J_m}),
\end{equation*}
where $\mu$ is the shear modulus and $J_m$ is a material stiffening parameter that imposes a maximum allowable stretch on the material. The first Piola-Kirchhoff stress is defined as
\begin{equation*}
P = -\frac{\mu J_m}{J_m-(I_1-2)}.
\end{equation*}

The Green-Lagrangian strains, $\epsilon_{ij} = \frac{1}{2}(F_{ki}F_{kj}-\delta_{ij})$, from finite element simulations were supplied to PINNs as reference data. All FEA analyses were carried out using FEniCS~\cite{Anders2010}.

\subsection{Sample complexity}

We used Delentropy~\cite{larkin2016, Khan2022} to analyze the spatial structure and pixel co-occurrence of the tissue image and evaluate the complexity of the tissue structural pattern. The Delentropy is expressed as 
\begin{equation*}
DE = -\frac{1}{2}\sum_{j=0}^{J-1}\sum_{i=0}^{I-1}p_{i,j}log_2p_{i.j},
\end{equation*} 
where $p_{i,j}$ is the deldensity probability density function and $I$ and $J$ are the number of discrete cells in the $x$ and $y$ dimensions of the probability density function. The deldensity probability density function is computed as 

\begin{equation*}
p_{i,j} = \frac{1}{4WH}\sum_{w=0}^{W-1}\sum_{h=0}^{H-1}\delta_{i, d_x}(w,h)\delta_{j, d_y}(w,h),
\end{equation*} 
where $W$ and $H$ are the width and height of the image, respectively. The Kronecker delta, $\delta$, describes the binning operation for generating a 2D histogram for the image. $d_x$ and $d_y$ are the derivative of the kernels in $x$ and $y$ directions, respectively. The local derivatives were estimated by computing the central difference corresponding to the convolution kernel of $[[-1,0,1],[-1,0,1],[-1,0,1]]$. 

\subsection{PINNs for heterogeneous material properties identification}
\label{pinns_framework}
As an extension to prior work on utilizing PINNs to determine the macroscopic constitutive parameters in linear and hyperelastic materials~\cite{haghighat2021, Wu2023}, we demonstrate a proof-of-concept study to identify the full-field elastic properties across the 2D geometric domain. Our work aims to facilitate future studies to better understand the interconnection between material microscopic properties and its macroscopic mechanical functions. All PINNs models were constructed using the DeepXDE library \cite{lu2021}, and the code will be made publicly available after publication in the GitHub repository \url{https://github.com/lu-group/pinn-heterogeneous-material}.

%\textcolor{blue}{\subsubsection{Optimization description}}

\subsubsection{Output transform}
\label{output_transform}
Rescale network outputs magnitudes to the scale of $\mathcal{O}(1)$ can help improve the optimization process. We normalized the displacement outputs by the mean and standard deviation obtained from the finite element analysis. Subsequently, we implemented hard constraint displacement boundary conditions in the same fashion described in Refs.~\cite{lu2021_2, Wu2023}. The transformation of the displacement vectors in a compact form is expressed as
\begin{equation*} 
u^*_i = \ell(\mathbf{x})(\mathcal{N}_{u^*_i} u_{i, \text{std}} +  u_{i, \text{mean}}), \\
\end{equation*}
where $u^*_i$ is the final PINN displacement vector, $\ell(\mathbf{x})$ is the hard constraint function, $\mathcal{N}_{u^*_i}$, is the normalized PINN displacement vector, $u_{i, \text{std}}$ and $u_{i, \text{mean}}$ are the mean and standard deviation of the reference displacement vectors. Further, we constrained the range of $E$ to physically realistic values using a sigmoid function
\begin{equation*} 
    E^* = \frac{4}{1 + e^{-\mathcal{N}_{E^*}}} + 1.
\end{equation*}

\subsubsection{Loss functions}
The training process aims to optimize the network parameters, $\theta_{\text{NN}}$, in
\begin{equation*}
\theta_{\text{NN}}^* = \underset{\theta_{\text{NN}}}{\arg\min} \mathcal{L} (\theta_{\text{NN}}),
\end{equation*}
where $\mathcal{L}(\theta_{\text{NN}})$ is the loss function that measures the total error concerning the PDEs (the momentum balance), constitutive relation, and reference data loss; the boundary condition loss is eliminated since we impose the boundary conditions using hard constraints. The total loss function for FCNN I and II is defined as
\begin{equation*}
\mathcal{L}(\theta_{\text{NN}}, \theta_{\text{mat}}) =
w_{\text{PDEs}}\mathcal{L}_{\text{PDEs}} +
w_{\text{Constitutive relations}}\mathcal{L}_{\text{Constitutive relations}} +
w_{\text{Data}}\mathcal{L}_{\text{Data}}, 
\end{equation*}
where $w_\bullet$ is the weight associated with its corresponding loss term $\mathcal{L}_\bullet$. The loss terms $\mathcal{L}_{\text{Data}}$ compute the mean squared error of predicted results on the collocation points and reference data, respectively. While displacement and strain data were available, we only used strains to facilitate data loss because including displacements did not improve model prediction. $\mathcal{L}_{\text{PDEs}}$ and $\mathcal{L}_{\text{Constitutive relations}}$ computes the mean squared error of the PDE residuals and constitutive behavior loss over the spatial domain. The weights associated with PDEs, Constitutive relation, and data loss were $1$, $E^2$, and $100$, respectively. The total loss function for FCNN III to V is reduced to 
\begin{equation*}
\mathcal{L}(\theta_{\text{NN}}, \theta_{\text{mat}}) =
w_{\text{PDEs}}\mathcal{L}_{\text{PDEs}} +
w_{\text{Data}}\mathcal{L}_{\text{Data}}.
\end{equation*}
Without loss of generality, we adopted the Swish activation function across all analyses and trained the network with 500,000 iterations with Adam optimizer and learning rate as 0.001 in all studies.  

\section*{Appendix A: Additional results on method comparison}

\begin{figure}[htbp]
\centering
\vspace{-2ex}
\includegraphics[width=1\textwidth]{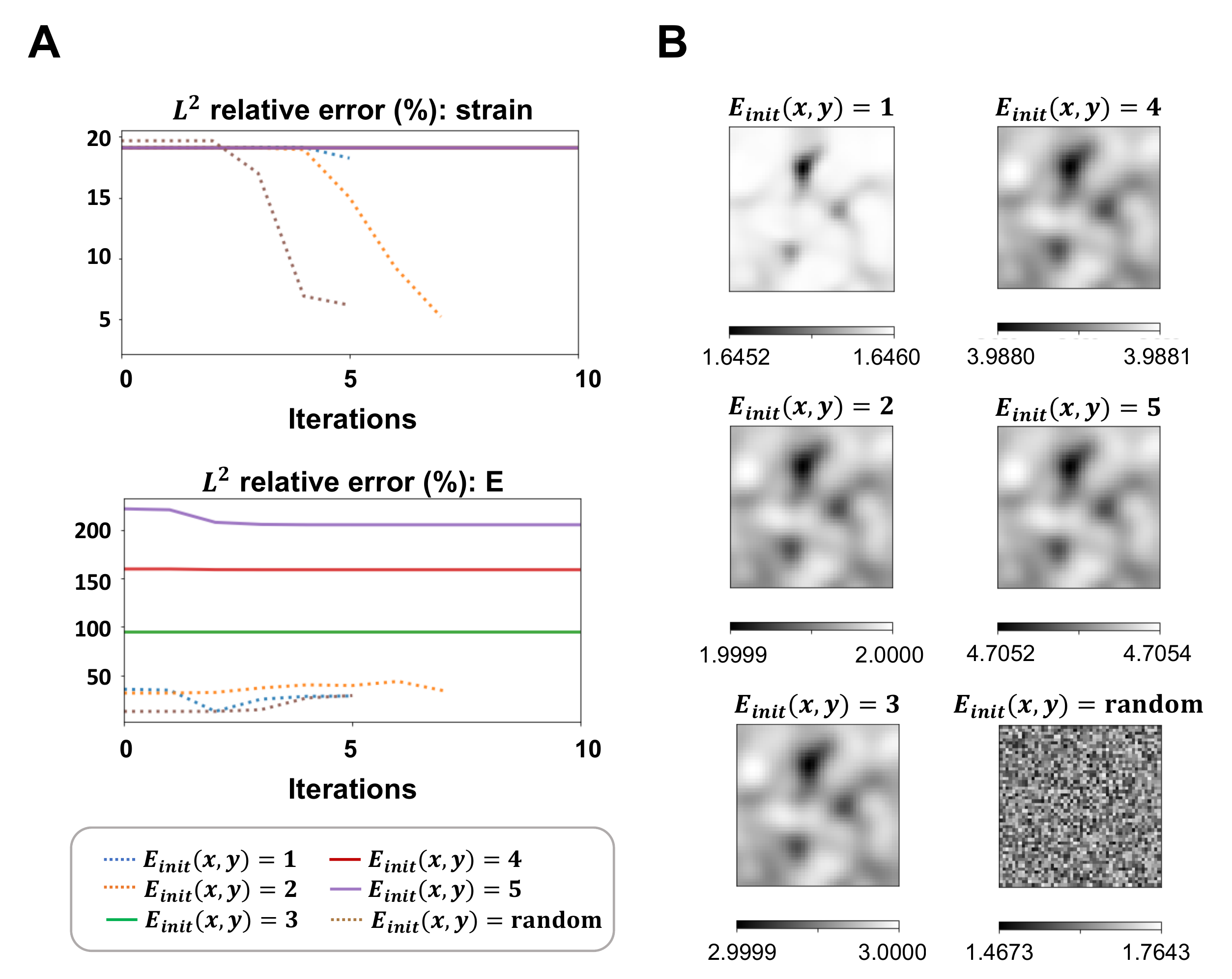}
\vspace{-2ex}
\caption{\textbf{Elasticity map convergence trajectories and best estimations.} (\textbf{A}) The $L^2$ relative errors convergence trajectories of strain and the elasticity map are shown. The $L^2$ relative errors in strain consistently remained around 20\% with initial values of $E_{\text{init}}(x, y) =$ 3, 4, and 5. The lack of convergence resulted in constant high errors in the elasticity maps, with $L^2$ relative errors exceeding 100\%. (\textbf{B}) The best estimates of the elasticity map from each elasticity map initialization are presented. The $L^2$ relative errors of the estimated elasticity map for $E_{\text{init}}(x, y)$ values of 1, 2, 3, 4, 5, and randomly distributed were 12.93\%, 31.84\%, 94.67\%, 158.87\%, 205.18\%, and 12.93\%.}\label{ifea_1}
\vspace{-2ex}
\end{figure}

We present inverse FEA results for the GRF example described in Section~\ref{re:setup}. The optimization process was facilitated using the adjoint method. Specifically, we set to identify a heterogeneous elasticity map that minimizes the $L^2$ relative errors of strain under the constraint that $1 \leq E \leq 5$. The maximum optimization iteration was set to 100, and the acceptable tolerance of the $L^2$ relative error for strains was set to $10^{-3}$. The convergence histories of strain and $E$ in the first 10 iterations are shown in Fig.~\ref{ifea_1}A.

We found that the initialization of the elasticity map, $E_{\text{init}}(x, y)$, has a significant influence on the estimation of the elasticity map. When $E_{\text{init}} =$ 3, 4, and 5, the $L^2$ relative errors of strain were consistently around 20\% throughout the optimization process. On the other hand, with $E_{\text{init}}(x, y) =$ 1, 2, and randomly distributed, the $L^2$ relative errors of strain decreased as the optimization progressed. However, the intermediate estimated elasticity field often resulted in divergence in the FEA, consequently leading to pre-mature termination of the optimization process.

When the initial elastic modulus $E_{\text{init}}(x, y)$ was set to 1, 2, and randomly distributed, a noticeable decrease in the $L^2$ relative errors of the strain was evident. However, this decline did not correspond to the $L^2$ relative errors of the elasticity map. Surprisingly, an increase in $L^2$ relative errors was observed in cases when $E_{\text{init}}(x, y)$ = 2 and randomly distributed, despite a decrease in strain error.

Fig.~\ref{ifea_1}B demonstrates the most accurate estimates, measured by the least $L^2$ relative errors, of elasticity maps resulting from the various $E_{\text{init}}(x, y)$. The results indicate that while the inverse FEA roughly uncovered the heterogeneous pattern, the range of elasticity moduli with $E_{\text{init}}(x, y)$ values of 1, 2, 3, 4, and 5 suggests a relatively uniform and close to homogeneous distribution. Notably, when the $E_{\text{init}}(x, y)$ values were randomly distributed, we observed a more reasonable range in elastic moduli. However, the resulting elasticity map estimation did not capture the true pattern of the elasticity map.

\begin{figure}[htbp]
\centering
\vspace{-2ex}
\includegraphics[width=1\textwidth]{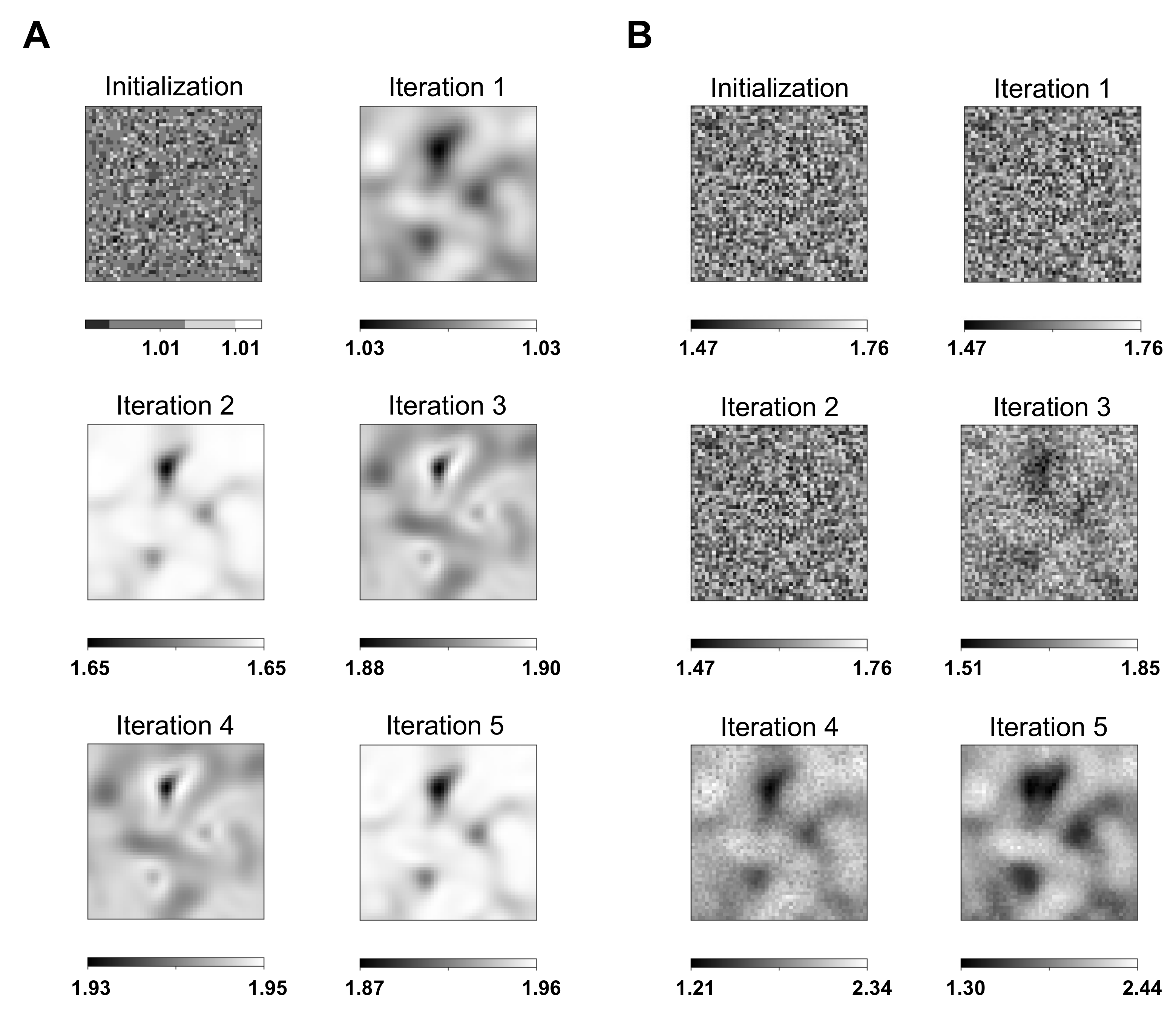}
\vspace{-2ex}
\caption{\textbf{Elasticity map evolutions.} The evolutions of the heterogeneous elasticity maps with $E_{\text{init}}(x, y) =$ 1 and randomly distributed as they provided the best estimates among the numerical experiments. (\textbf{A}) The $L^2$ relative errors of the elasticity map at each iteration, when $E_{\text{init}}(x, y) =$ 1, were as follows: 35.75\%, 34.80\%, 12.93\%, 25.49\%, 28.45\%, and 28.79\%. (\textbf{B}) The $L^2$ relative errors of the elasticity map at each iteration, with randomly distributed $E_{\text{init}}(x, y)$, were 12.93\%, 12.93\%, 12.93\%, 14.69\%, 26.78\%, and 29.22\%.}\label{ifea_2}
\vspace{-2ex}
\end{figure}

Fig.~\ref{ifea_2} presents the evolutions of the heterogeneous elasticity maps with $E_{\text{init}}(x, y) =$ 1 and randomly distributed. We observed that the optimizer struggled to uncover the true range of the elastic moduli in the case when initializing the elasticity map uniformly to 1. On the contrary, the random initialization of the elasticity map seemed to better account for the wide range of elastic moduli.

\clearpage
\section*{Appendix B: Additional results on the selected architectures}
\begin{figure}[htbp]
\centering
\vspace{-2ex}
\includegraphics[width=1\textwidth]{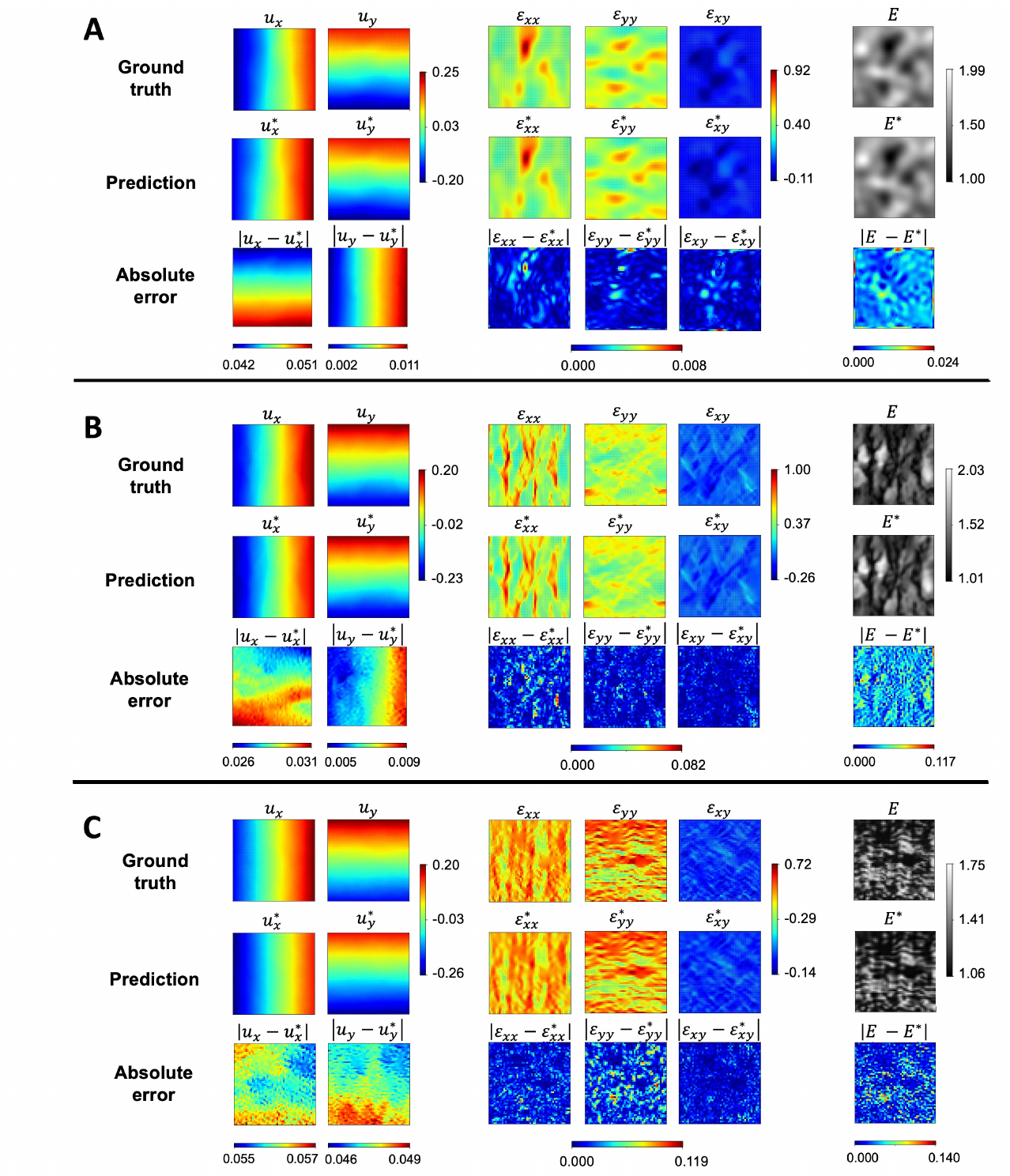}
\vspace{-2ex}
\caption{\textbf{PINN prediction and pointwise error.} The FEA ground truth, PINN estimation (denoted with superscript $*$), and the absolute errors of displacement, strain, and elastic modulus distribution for examples (\textbf{A}), (\textbf{B}), and (\textbf{C}) are shown. The PINN architecture effectively resolves the local features of complex strain and elastic modulus fields. The absolute errors were up to $\mathcal{O}(10^{-2})$ for the displacement and up to $\mathcal{O}(10^{-1})$ for the strain and elasticity fields.}\label{best}
\vspace{-2ex}
\end{figure}

Fig.~\ref{best} presents the estimated displacements, strain, and elastic modulus distribution of the GRF, brain tissue, and tricuspid valve tissues examples using network architecture IIB. The absolute relative errors between the PINN estimation and FEA ground truth are also provided for comparison. Despite the displacement boundaries being unconstrained in this network architecture, PINN model IIB demonstrated highly accurate estimation of the displacement, strain, and elastic parameter distribution across the three examples with relative errors up to $\mathcal{O}(10^{-2})$ for full-field displacements and strains and $\mathcal{O}(10^{-1})$ for the elastic moduli. These results indicate the consistency and generality of the proposed neural network architecture for estimating the mechanical properties of complex microstructures.

As the elasticity fields become more complex, with increasing local features, the corresponding strain fields become more complex. Consequently, we observed higher absolute errors associated with both the strain and elastic modulus in the second and third examples. In the first example, the maximum absolute errors of the strain and elasticity field were 0.01 and 0.02, respectively. However, these errors were 0.08 and 0.12 in the second example, and 0.12 and 0.14 in the third example. It is worth noting that increasing the size of the FCNN architecture may help resolve the finer features in the elasticity map, leading to better accuracy in the more complex examples.

\clearpage
\section*{Acknowledgments}
%TC:ignore
This work was supported by the Cora Topolewski Pediatric Valve Center at the Children's Hospital of Philadelphia, an Additional Ventures Expansion Award, the National Institutes of Health, and the U.S. Department of Energy. WW was initially supported by NHLBI T32 HL007915 and transitioned to NHLBI K25 HL168235. MAJ was supported by NIH R01 HL153166, an Additional Ventures Single Ventricle Grant, and the Topolewski Endowed Chair in Pediatric Cardiology. MD and LL were supported by the U.S. Department of Energy [DE-SC0022953]. 
%TC:endignore

\section*{Competing interests}
All authors declare no financial or non-financial competing interests. 

\section*{Code availability}
The code for this study is available after publication in the GitHub repository \url{https://github.com/lu-group/pinn-heterogeneous-material}.

\bibliography{library}

\end{document}